\documentclass[12pt,thmsa]{article}
\usepackage{amsmath}
\usepackage{amssymb}
\usepackage{sw20lart}

\input{tcilatex}
\begin{document}

\author{Kahar El-Hussein \\
Department of Mathematics, Faculty of Science and Arts at Al Qurayat,\\
Al Jouf University, KSA \& Forat University, Deir el Zore, Syria\\
E-mail: kumath@ju.edu.sa, kumath@hotmail.com}
\title{Plancherel Theorem on the Symplectic Group $SP(4,\mathbb{R})$}
\maketitle

\begin{abstract}
Let $SL(4,\mathbb{R})$ be the $15-$ dimensional connected semisimple Lie
group and let $SL(4,\mathbb{R})=KAN$ \ be the Iwasawa decomposition. \ Let $%
\mathbb{R}^{4}$ $\rtimes SL(4,\mathbb{R})$ be the group of the semidirect
product of $SL(4,\mathbb{R})$ with the real $\ $vector group $\mathbb{R}^{4}.
$ The goal of this paper is to define the Fourier transform on $SL(4,\mathbb{%
R})$\ in order to obtain the Plancherel theorem on $SL(4,\mathbb{R})$ and so
on $\mathbb{R}^{4}$ $\rtimes SL(4,\mathbb{R}).$ Since the symplectic group $%
SP(4,\mathbb{R})$ is a subgroup of $SL(4,\mathbb{R})$, then it will be easy
to get the Plancherel theorem on $SP(4,\mathbb{R})$ and so on its
inhomogeneous group. To this end, we obtain some interesting results on its
nilpotent symplectic group
\end{abstract}



\section{\protect\bigskip Introduction}

\textbf{1. }As well known the connected semisimple Lie group $SL(n,\mathbb{R}%
)$, consists of the following matrices 
\begin{equation}
SL(n,\mathbb{R})=\{A=GL(n,\mathbb{R)};\det A=1\}
\end{equation}

The\bigskip\ group $Sp(2n,\mathbb{R})$ is a subgroup of $SL(2n,\mathbb{R}),$
which is%
\begin{equation}
SP(2n,\mathbb{R})=\{g\in SL(2n,\mathbb{R});g\jmath g^{t}=\jmath \}
\end{equation}%
where $\jmath $ is the symplectic matrix defined by

\begin{equation}
\jmath =\left( 
\begin{array}{cc}
0 & I_{n} \\ 
-I_{n} & 0%
\end{array}%
\right)
\end{equation}%
and $0$ and $I$ are the $n\times n$ zero and identity matrices. It is clear $%
det$ $\jmath =1,$ $\jmath ^{2}=I,$ and $\jmath ^{t}=\jmath ^{-1}=-\jmath $

\bigskip The $10-$demesional symplectic group $SP(4,\mathbb{R)}$. If $g$ $%
\in SP(4,\mathbb{R)}$, then%
\begin{equation}
g=\left( 
\begin{array}{cccc}
x_{11} & x_{12} & x_{13} & x_{14} \\ 
x_{21} & x_{22} & x_{14} & x_{24} \\ 
x_{31} & x_{32} & -x_{11} & -x_{21} \\ 
x_{32} & x_{42} & -x_{12} & -x_{22}%
\end{array}%
\right) ,x_{ij}\in \mathbb{R}\text{, }1\leq i,j\leq 4
\end{equation}

The goal of this paper, is to define the Fourier transform in order to
obtain the Plancherel theorem on the group $SP(4,\mathbb{R)}$ and its
inhomogeneous group. Therefore, I will define the Fourier transform on $SL(4,%
\mathbb{R)}$, and I will prove its Plancherel theorem. Besides, I will
demonstrate the existence theorem and hypoellipticity of the partial
differential equations on its nilpotent symplectic

\section{\protect\bigskip Notation and Results}

\textbf{2.1. }In the following and far away from the representations theory
of Lie groups we use the Iwasawa decomposition of $SL(4,\mathbb{R}),$ to
define the Fourier transform and to get the Plancherel formula on the
connected real semisimple Lie group $SL(4,\mathbb{R}).$ Therefore let $SL(4,%
\mathbb{R})$\ be the complex Lie group, which is

\begin{equation}
SL(4,\mathbb{R})=\{\left( A=%
\begin{array}{cccc}
a_{11} & a_{12} & a_{13} & a_{14} \\ 
a_{21} & a_{22} & a_{23} & a_{24} \\ 
a_{31} & a_{32} & a_{33} & a_{34} \\ 
a_{41} & a_{42} & a_{43} & a_{44}%
\end{array}%
\right) :\text{ }a_{ij}\in \mathbb{R},1\leq i,j\leq 4\text{ \textit{and} }%
\det A=1\}
\end{equation}%
Let $G=SL(4,\mathbb{R})=KNA$ be the Iwasawa decomposition of $G$, where 
\begin{eqnarray}
K &=&SO(4)  \notag \\
N &=&\{\text{ }\left( 
\begin{array}{cccc}
1 & x_{1} & x_{3} & x_{6} \\ 
0 & 1 & x_{2} & x_{5} \\ 
0 & 0 & 1 & x_{4} \\ 
0 & 0 & 0 & 1%
\end{array}%
\right) _{i}x_{i}\in \mathbb{R}\text{ },1\leq i\leq 6\}  \notag \\
A &=&\{\left( 
\begin{array}{cccc}
a_{1} & 0 & 0 & 0 \\ 
0 & a_{1} & 0 & 0 \\ 
0 & 0 & a_{1} & 0 \\ 
0 & 0 & 0 & a_{1}%
\end{array}%
\right) :\text{ }a_{i}\in \mathbb{R}_{+}^{\star },1\leq i\leq
4,a_{1}a_{2}a_{3}a_{4}=1\}
\end{eqnarray}

\bigskip Hence every $g\in SL(4,\mathbb{R})$ can be written as $g=kan\in
SL(4,\mathbb{R}),$ where $k\in K,$ $a\in A,$ $n\in N.$ We denote by $%
L^{1}(SL(4,\mathbb{R}))$ the Banach algebra that consists of all complex
valued functions on the group $SL(4,\mathbb{R})$, which are integrable with
respect to the Haar measure $dg$ of $SL(4,\mathbb{R})$ and multiplication is
defined by convolution product on $SL(4,\mathbb{R})$ , and we denote by $%
L^{2}(SL(4,\mathbb{R}))$ the Hilbert space of $SL(4,\mathbb{R})$. So we have
for any $f\in L^{1}(SL(4,\mathbb{R}))$ and $\phi \in L^{1}(SL(4,\mathbb{R}))$
\begin{equation}
\phi \ast f(h)=\int\limits_{G}f(g^{-1}h)\phi (g)dg
\end{equation}%
The Haar measure $dg$ on a connected real semi-simple Lie group $G$ $=SL(4,%
\mathbb{R)}$, can be calculated from the Haar measures $dn,$ $da$ and $dk$
on $N;A$ and $K;$respectively, by the formula%
\begin{equation}
\int\limits_{SL(4,\mathbb{R})}f(g)dg=\int\limits_{A}\int\limits_{N}\int%
\limits_{K}f(ank)dadndk
\end{equation}

\textbf{2.2.} Keeping in mind that $a^{-2\rho }$ is the modulus of the
automorphism $n\rightarrow $ $ana^{-1}$ of $N$ we get also the following
representation of $dg$ 
\begin{equation}
\int\limits_{SL(4,\mathbb{R})}f(g)dg=\int\limits_{A}\int\limits_{N}\int%
\limits_{K}f(ank)dadndk=\int\limits_{N}\int\limits_{A}\int%
\limits_{K}f(nak)a^{-2\rho }dndadk
\end{equation}%
where%
\begin{equation*}
\rho =2^{-1}\sum_{\alpha \rangle 0}m(\alpha )\alpha
\end{equation*}%
and $m(\alpha )$ denotes the multiplicity of the root $\alpha $ and $\rho =$
the dimension of the nilpotent group $N.$ Furthermore, using the relation $%
\int\limits_{G}f(g)dg=\int\limits_{G}f(g^{-1})dg,$ we receive 
\begin{equation}
\int\limits_{SL(4,\mathbb{R})}f(g)dg=\int\limits_{K}\int\limits_{A}\int%
\limits_{N}f(kan)a^{2\rho }dndadk
\end{equation}

\section{\protect\bigskip \textbf{Fourier Transform and Plancherel Theorem On%
} $N$}

\textbf{3.1.} Let $N$ be the real group consisting of all matrices of the
form 
\begin{equation}
\left( 
\begin{array}{cccc}
1 & x_{1} & x_{3} & x_{6} \\ 
0 & 1 & x_{2} & x_{5} \\ 
0 & 0 & 1 & x_{4} \\ 
0 & 0 & 0 & 1%
\end{array}%
\right)
\end{equation}%
where $(x_{1},x_{2},$ $x_{3},$ $x_{4},$ $x_{5},x_{6})\in \mathbb{R}^{6}$ .
The group can be identified with the group $(\mathbb{R}^{3}\underset{\rho
_{2}}{\rtimes }\mathbb{R}^{2}$ $)\underset{\rho _{1}}{\rtimes }\mathbb{R}$\
be the semidirect product of the real vector groups $\mathbb{R}$, $\mathbb{R}%
^{2}$ and $\mathbb{R}^{3}$, where $\rho _{2}$ is the group homomorphism 
\hspace{0.05in}$\rho _{2}:\mathbb{R}^{2}\rightarrow Aut(\mathbb{R}^{3}),$
which is defined by 
\begin{equation}
\rho _{2}(x_{3},x_{2}\text{ }%
)(y_{6},y_{5},y_{4})=(y_{6}+x_{3}y_{4},y_{5}+x_{2}y_{4},y_{4})
\end{equation}%
and $\rho _{1}$ is the group homomorphism\hspace{0.05in}$\rho _{1}:\mathbb{R}%
\rightarrow Aut(\mathbb{R}^{3}\underset{\rho _{2}}{\rtimes }\mathbb{R}^{2}),$
which is given by 
\begin{equation}
\rho _{1}(x_{1}\text{ }%
)(y_{6},y_{5},y_{4},y_{3},y_{2})=(y_{6}+x_{1}y_{5},y_{5},y_{4},y_{3}+x_{1}y_{2},y_{2})
\end{equation}%
where $Aut(\mathbb{R}^{3})$ $(resp.Aut(\mathbb{R}^{3}\underset{\rho _{2}}{%
\rtimes }\mathbb{R}^{2}))$ is the group of all automorphisms of $\mathbb{R}%
^{3}$ $(resp.(\mathbb{R}^{3}\underset{\rho _{2}}{\ltimes }\mathbb{R}^{2})),$
see $[6].$

Let $L=\mathbb{R}^{3}\times \mathbb{R}^{2}\times \mathbb{R}^{2}\mathbb{%
\times R}\times \mathbb{R}$\ be the group with law: 
\begin{eqnarray}
&&X.Y=(x_{6},x_{5},x_{4},x_{3},x_{2},t_{3},t_{2},x_{1},t_{1})(y_{6},y_{5},y_{4},y_{3},y_{2},s_{3},s_{2},y_{1},s_{1})
\notag \\
&=&((x_{6},x_{5},x_{4},x_{3},x_{2},t_{3},t_{2})(\rho
_{1}(t_{1})(y_{6},y_{5},y_{4},y_{3},y_{2},s_{3},s_{2}),y_{1}+x_{1},s_{1}+t_{1})
\notag \\
&=&((x_{6},x_{5},x_{4},x_{3})\rho
_{2}(t_{3},t_{2})(y_{6}+t_{1}y_{5},y_{5},y_{4},y_{3},s_{3}+t_{1}s_{2},s_{2}),x_{2}+y_{2},y_{1}+x_{1},s_{1}+t_{1})
\notag \\
&=&((x_{6},x_{5},x_{4})+(y_{6}+t_{1}y_{5}+t_{3}y_{4},y_{5}+t_{2}y_{4},y_{4}),x_{3}+y_{3},s_{3}+t_{1}s_{2},
\notag \\
&&s_{2}+t_{2},x_{2}+y_{2},y_{1}+x_{1},s_{1}+t_{1})  \notag \\
&=&(x_{6}+y_{6}+t_{1}y_{5}+t_{3}y_{4},x_{5}+y_{5}+t_{2}y_{4},x_{4}+y_{4},x_{3}+y_{3},t_{3}+s_{3}+t_{1}s_{2},
\notag \\
&&y_{2}+x_{2},s_{2}+t_{2},y_{1}+x_{1},s_{1}+t_{1})
\end{eqnarray}%
for all $(X,Y)\in L^{2}.$ In this case the group $N$ can be identified with
the closed subgroup $\mathbb{R}^{3}\times \{0\}\underset{\rho _{1}}{\rtimes }%
\mathbb{R}^{2}\mathbb{\times }$ $\{0\}\ \underset{\rho _{1}}{\rtimes }%
\mathbb{R}$ of $L$ and $B$ with the closed subgroup $\mathbb{R}^{3}\times 
\mathbb{R}^{2}\times \{0\}\mathbb{\times R}\times \{0\}$\ of $L,$ where $B=%
\mathbb{R}^{3}\times \mathbb{R}^{2}\mathbb{\times }$ $\mathbb{R}$ the group,
which is the direct product of the real vector groups $\mathbb{R}$, $\mathbb{%
R}^{2}$ and $\mathbb{R}^{3}$

Let $C^{\infty }(N),$ $\mathcal{D}(N),$ $\mathcal{D}^{\prime }(N),$ $%
\mathcal{E}^{\prime }(N)$ be the space of $C^{\infty }$- functions, $%
C^{\infty }$ with compact support, distributions and distributions with
compact support on $N$ respectively$.$We denote by $L^{1}(N)$ the Banach
algebra that consists of all complex valued functions on the group $N$,
which are integrable with respect to the Haar measure of $N$ and
multiplication is defined by convolution on $N$, and we denote by$L^{2}(N)$
the Hilbert space of $N$. \newline

\textbf{Definition 3.1.}\textit{\ For every }$f\in C^{\infty }(N)$, \textit{%
one can define function} $\widetilde{f}\in C^{\infty }(L)$ \textit{as
follows:} 
\begin{equation}
\widetilde{f}(x,x_{3},x_{2},t_{3},t_{2},x_{1},t_{1})=f(\rho _{1}(x_{1})(\rho
_{2}(x_{3},x_{2})(x),t_{3}+x_{3},t_{2}+x_{2}),t_{1})
\end{equation}%
\textit{for all} $(x,x_{3},x_{2},t_{3},t_{2},x_{1},t_{1})\in L,$ here $%
x=(x_{6},x_{5},x_{4})\in \mathbb{R}^{3}.$

\textbf{Remark 3.1. }\textit{The function} $\widetilde{f}$ \textit{is
invariant in the following sense:} 
\begin{eqnarray}
&&\widetilde{f}((\rho _{1}(h)((\rho
_{2}(r,k)(x),x_{3}-r,x_{2}-k,t_{3}+r,t_{2}+k),x_{1}-h,t_{1}+h)  \notag \\
&=&\widetilde{f}(x,x_{3},x_{2},t_{3},t_{2},x_{1},t_{1})
\end{eqnarray}%
\textit{for any} $(x,x_{3},x_{2},t_{3},t_{2},x_{1},t_{1})\in L$\textit{, }$%
h\in \mathbb{R}$ \textit{and }$(r,k)\in \mathbb{R}$, where $%
x=(x_{6},x_{5},x_{4})\in \mathbb{R}^{3}$. \textit{So every function} $\psi
(x,x_{3},x_{2},x_{1})$ \textit{on} $N$\textit{\ extends uniquely as an
invariant function} $\widetilde{\psi }%
(x,x_{3},x_{2},t_{3},t_{2},x_{1},t_{1}) $ \textit{on} $L.$

\textbf{Theorem 3.1.}\textit{\ For every function} $F\in C^{\infty }(L)$ 
\textit{invariant in sense} $(16)$ \textit{and for every }$\varphi \mathcal{%
\in }$ $\mathcal{D}(N)$, \textit{we have} 
\begin{equation}
u\text{ }\ast F(x,x_{3},x_{2},t_{3},t_{2},x_{1},t_{1})=u\text{ }\ast
_{c}F(x,x_{3},x_{2},t_{3},t_{2},x_{1},t_{1})
\end{equation}%
\textit{for every} $(X,x_{3},x_{2},t_{3},t_{2},x_{1},t_{1})$ $\in L$, 
\textit{where} $\ast $\textit{\ signifies the convolution product on} $N$\ 
\textit{with respect the variables} $(x,t_{3},t_{2},t_{1})$\ \textit{and} $%
\ast _{c}$\textit{signifies the commutative convolution product on} $B$ 
\textit{with respect the variables }$(x,x_{3},x_{2},x_{1}).$

\textit{Proof}\textbf{:}\textit{\ }\textbf{\ }In fact we have\ 
\begin{eqnarray}
&&\varphi \ast F(x,x_{3},x_{2},t_{3},t_{2},x_{1},t_{1})  \notag \\
&=&\int\limits_{N}F\left[
(y,y_{3},y_{2},s)^{-1}(X,x_{3},x_{2},t_{3},t_{2},x_{1},t_{1})\right]
u(y,y_{3},y_{2},s)dydy_{3}dy_{2}ds  \notag \\
&=&\int\limits_{N}F\left[ (\rho
_{1}(s^{-1})(y,y_{3},y_{2})^{-1},-s)(x,x_{3},x_{2},t_{3},t_{2},x_{1},t_{1})%
\right] u(y,y_{3},y_{2},s)dydy_{3}dy_{2}ds  \notag \\
&=&\int\limits_{N}F[(\rho _{1}(s^{-1})((\rho
_{2}(y_{3},y_{2})^{-1}((-y)+(x))),x_{3},x_{2},t_{3}-y_{3},t_{2}-y_{2}),x_{1},t_{1}-s)]
\notag \\
&&u(y,y_{3},y_{2},s)dydy_{3}dy_{2}ds
\end{eqnarray}%
Since $F$ is invariant in sense $(16),$ then for every $%
(x,x_{3},x_{2},t_{3},t_{2},x_{1},t_{1})\in L$ we get 
\begin{eqnarray}
&&\varphi \ast F(x,x_{3},x_{2},t_{3},t_{2},x_{1},t_{1})  \notag \\
&=&\int\limits_{N}F[(\rho _{1}(s^{-1})(\rho
_{2}(y_{3},y_{2})^{-1}(-y+x),x_{3},x_{2},t_{3}-y_{3},t_{2}-y_{2}),  \notag \\
&&x_{1},t_{1}-s)]u(y,y_{3},y_{2},s)dydy_{3}dy_{2}ds  \notag \\
&=&\int\limits_{N}F\left[ x-y,x_{3}-y_{3},x_{2}-y_{2},t_{3},t_{2},x_{1}-s,t%
\right] u(y,y_{3},y_{2},s)dydy_{3}dy_{2}ds  \notag \\
&=&\varphi \ast _{c}F(x,x_{3},x_{2},t_{3},t_{2},x_{1},t_{1})
\end{eqnarray}%
\textbf{\ \ \ \ \ \ \ \ \ \ \ \ \ \ \ \ \ \ \ \ \ \ \ \ \ \ \ \ \ \ \ \ \ \
\ \ \ \ \ \ \ \ }

\textbf{\ \ \ \ \ \ \ \ \ \ }

As in $[6]$, we will define the Fourier transform on $G$. Therefore let $%
\mathcal{S}(N)$ be the Schwartz space of $N$ which\hspace{0.05in}can be
considered as the Schwartz space of $\mathcal{S}(B),$ and let $\mathcal{S}%
^{\prime }(N)$ be the space of all tempered distributions on $N.$

\textbf{Definition 3.2.}\textit{\ If }$f\in $\textit{\ }$\mathcal{S}(N),$%
\textit{\ one can define its Fourier transform }$\mathcal{F}f$\textit{\ 
\hspace{0.05in}by the Fourier transform on its vector group}: 
\begin{equation}
\mathcal{F}f\text{ \ }(\xi )=\int\limits_{N}f(X)\text{ }e^{-\text{ }i\text{ }%
\langle \xi \text{,}X\rangle }\text{ }dX
\end{equation}%
\textit{for any} $\xi =(\xi _{6},\xi _{5},\xi _{4},\xi _{3},\xi _{2},\xi
_{1})\in \mathbb{R}^{6},$ and $X=(x_{6},x_{5},x_{4},x_{3},x_{2},x_{1})\in 
\mathbb{R}^{6}$, \textit{where} $\langle \xi $,$X\rangle =\xi _{6}x_{6}+\xi
_{5}x_{5}+\xi _{4}x_{4}+\xi _{3}x_{3}+\xi _{2}x_{2}+\xi _{1}x_{1}$ \textit{%
and} $dX=dx_{6}dx_{5}dx_{4}dx_{3}dx_{2}dx_{1}$\textit{is the Haar measure on}
$N$. \textit{The mapping}\ $f$ \ $\rightarrow $ \ $\mathcal{F}f$ \ \textit{%
is isomorphism of the topological vector space} \ $\mathcal{S}(N)$ \textit{%
onto} $\mathcal{S}(\mathbb{R}^{6}).$

\textbf{Theorem 3.2. }\textit{The Fourier transform }$\mathcal{F}$ \textit{%
satisfies :} 
\begin{equation}
\overset{\vee }{\varphi }\ast f(0)=\int\limits_{\mathbb{R}^{4}}\mathcal{F}%
f(\xi )\overline{\mathcal{F}u(\xi )}d\xi
\end{equation}%
\textit{for every }$f\in S(N)$\textit{\ and }$\varphi \in S(N),$ \textit{%
where\ }$\overset{\vee }{\varphi }(X)=\overline{u(X^{-1})},$\textit{\hspace{%
0.05in}}$\xi =(\xi _{6},\xi _{5},\xi _{4},\xi _{3},\xi _{2},\xi _{1}),$%
\textit{\ }$d\xi =d\xi _{1}d\xi _{2}d\xi _{3}d\xi _{4}d\xi _{5}d\xi _{6},$%
\textit{\ is the Lebesgue\ measure on }$\mathbb{R}^{6},$\textit{\ and } $%
\ast $\textit{\ denotes the convolution product on }$N$\textit{.} \medskip

\textit{Proof:\ }By the classical Fourier transform, we have: 
\begin{eqnarray}
&&\overset{\vee }{\varphi }\ast f(0)=\int\limits_{\mathbb{R}^{4}}\mathcal{F}(%
\overset{\vee }{\varphi }\ast f)(\xi )d\xi  \notag \\
&=&\int\limits_{\mathbb{R}^{6}}\int\limits_{N}\overset{\vee }{\varphi }\ast
f(X)\text{ }e^{-i\left\langle \xi ,X\right\rangle }\text{ }dXd\xi  \notag \\
&=&\int\limits_{\mathbb{R}^{6}}\int\limits_{N}\int\limits_{N}f(YX)\overline{%
u(Y)}e^{-i\left\langle \xi ,X\right\rangle }\text{ }dY\text{ }dXd\xi .
\end{eqnarray}

By change of variable $YX=X^{\prime }$ with $%
Y=(x_{6},x_{5},x_{4},x_{3},x_{2},x_{1})$ and $X^{\prime
}=(y_{6},y_{5},y_{4},y_{3},y_{2},y_{1})$, we get 
\begin{eqnarray*}
X &=&Y^{-1}X^{\prime
}=(x_{6},x_{5},x_{4},x_{3},x_{2},x_{1})^{-1}(y_{6},y_{5},y_{4},y_{3},y_{2},y_{1})
\\
&=&(y_{6}-x_{6}+x_{1}x_{5}-x_{1}y_{5}-x_{1}x_{2}x_{4}+x_{3}x_{4}-x_{3}y_{4}+x_{1}x_{2}y_{4},y_{5}-x_{5}+x_{2}y_{4}-x_{2}x_{4},
\\
&&x_{4}+y_{4},y_{3}-x_{3}-x_{1}y_{2}+x_{1}x_{2},y_{2}-x_{2},y_{1}-x_{1})
\end{eqnarray*}%
and 
\begin{eqnarray*}
&&-i\left\langle \xi ,X\right\rangle \\
&=&-i\left\langle \xi ,Y^{-1}X^{\prime }\right\rangle \\
&=&-i[(y_{6}-x_{6}+x_{1}x_{5}-x_{1}y_{5}-x_{1}x_{2}x_{4}+x_{3}x_{4}-x_{3}y_{4}+x_{1}x_{2}y_{4})\xi _{6}+(y_{5}-x_{5}+x_{2}y_{4})\xi _{5}
\\
&&-x_{2}x_{4}\xi _{5}+(y_{4}-x_{4})\xi
_{4}+(y_{3}-x_{3}-x_{1}y_{2}+x_{1}x_{2})\xi _{3}+(y_{2}-x_{2})\xi
_{2}+(y_{1}-x_{1})\xi _{1}] \\
&=&-i[(x_{6}\xi _{6}-y_{6}\xi _{6})+(-x_{2}x_{4}\xi _{6}+x_{2}y_{4}\xi
_{6}-y_{2}\xi _{3}+x_{2}\xi _{3}+x_{5}\xi _{6}-y_{5}\xi _{6}-\xi
_{1})x_{1}-y_{1}\xi _{1} \\
&&+(y_{5}\xi _{5}-x_{5}\xi _{5})+(y_{4}\xi _{5}-x_{4}\xi _{5}-\xi
_{2})x_{2}+y_{2}\xi _{2}+y_{3}\xi _{3}+(x_{4}\xi _{6}-y_{4}\xi _{6}-\xi
_{3})x_{3}+(y_{4}-x_{4})\xi _{4}
\end{eqnarray*}

So we obtain

\begin{eqnarray*}
&&e^{-i(y_{6}-x_{6}+x_{1}x_{5}-x_{1}y_{5}-x_{1}x_{2}x_{4}+x_{3}x_{4}-x_{3}y_{4}+x_{1}x_{2}y_{4})\xi _{6}}%
\text{ }e^{-i((y_{5}-x_{5}+x_{2}y_{4}-x_{2}x_{4})\xi _{5}+(y_{4}-x_{4})\xi
_{4})} \\
&&e^{-i((y_{3}-x_{3}-x_{1}y_{2}+x_{1}x_{2})\xi _{3}+(y_{2}-x_{2})\xi
_{2}+(y_{1}-x_{1})\xi _{1})}\text{ } \\
&=&e^{-i((x_{6}\xi _{6}-y_{6}\xi _{6})+(-x_{2}x_{4}\xi _{6}+x_{2}y_{4}\xi
_{6}-y_{2}\xi _{3}+x_{2}\xi _{3}+x_{5}\xi _{6}-y_{5}\xi _{6}-\xi
_{1})x_{1}-y_{1}\xi _{1})}\text{ } \\
&&e^{-i((y_{5}\xi _{5}-x_{5}\xi _{5})+(y_{4}\xi _{5}-x_{4}\xi _{5}-\xi
_{2})x_{2}+y_{2}\xi _{2})-i(y_{3}\xi _{3}+(x_{4}\xi _{6}-y_{4}\xi _{6}-\xi
_{3})x_{3}+(y_{4}-x_{4})\xi _{4})}\text{ }
\end{eqnarray*}

By the invariance of the Lebesgue,s measures $d\xi _{1},d\xi _{2}$ \hspace{%
0.05in}and $d\xi _{3}$, we get 
\begin{eqnarray*}
&&\overset{\vee }{\varphi }\ast
f(0)=\int\limits_{N}\int\limits_{N}\int\limits_{\mathbb{R}^{6}}f(X^{\prime
})e^{-i((x_{6}\xi _{6}-y_{6}\xi _{6})+(-x_{2}x_{4}\xi _{6}+x_{2}y_{4}\xi
_{6}-y_{2}\xi _{3}+x_{2}\xi _{3}+x_{5}\xi _{6}-y_{5}\xi _{6}-\xi
_{1})x_{1}+y_{1}\xi _{1})} \\
&&e^{-i((y_{5}\xi _{5}-x_{5}\xi _{5})+(y_{4}\xi _{5}-x_{4}\xi _{5}-\xi
_{2})x_{2}+y_{2}\xi _{2})}e^{-i(y_{3}\xi _{3}+(x_{4}\xi _{6}-y_{4}\xi
_{6}-\xi _{3})x_{3}+(y_{4}-x_{4})\xi _{4})}\overline{\varphi (Y)}dY\hspace{%
0.05in}dX^{\prime }d\xi \\
&=&\int\limits_{N}\int\limits_{N}\int\limits_{\mathbb{R}^{6}}f(X^{\prime
})e^{-i(y_{6}\xi _{6}-x_{6}\xi _{6}+y_{5}\xi _{5}-x_{5}\xi _{5}+y_{4}\xi
_{4}-x_{4}\xi _{4}+y_{3}\xi _{3}-\xi _{3}x_{3}-\xi _{2}x_{2}+y_{2}\xi
_{2}-\xi _{1}x_{1}+y_{1}\xi _{1})}\overline{\varphi (Y)}dY\hspace{0.05in}%
dX^{\prime }d\xi \\
&=&\int\limits_{\mathbb{R}^{6}}\mathcal{F}f(\xi )\text{ }\overline{\mathcal{F%
}\varphi (\xi )}d\xi
\end{eqnarray*}%
for any $Y$ $=(x_{6},x_{5},x_{4},x_{3},x_{2},x_{1})\in $ $\mathbb{R}^{6}$
and $X^{\prime }=(y_{6},y_{5},y_{4},y_{3},y_{2},y_{1})\in \mathbb{R}^{6},$
where $0=(0,0,0,0,0,0)$ is the identity of $N$. The theorem is proved.

\textbf{Corollary 3.1. }\textit{In theorem \textbf{3.2.} if\hspace{0.05in}we
replace }$\varphi $ \textit{by }$f,$\textit{\ we obtain the Plancherel,s
formula on }$N$ 
\begin{equation}
\overset{\vee }{f}\ast f(0)=\int\limits_{N}\left\vert f(X)\right\vert
^{2}dX=\int\limits_{\mathbb{R}^{6}}\left\vert \mathcal{F}f(\xi )\right\vert
^{2}d\xi
\end{equation}

\section{\textbf{Fourier Transform and Plancherel Theorem On} $SL(4,\mathbb{%
R)}$}

\textbf{4.1. }Let $\underline{k}$ be the Lie algebra of $K=SO(4)$. Let $%
(X_{1},X_{2},X_{3},X_{4})$ a basis of $\underline{k}$ , such that the both
operators%
\begin{equation}
\Delta =\sum\limits_{i=1}^{3}X_{i}^{2},D_{q}=\sum\limits_{0\leq l\leq
q}\left( -\sum\limits_{i=1}^{3}X_{i}^{2}\right) ^{l}
\end{equation}%
are left and right invariant (bi-invariant) on $K,$ this basis exist see $%
[4, $ $p.564)$. For $l\in 
\mathbb{N}
$. Let $D^{l}=(1-\Delta )^{l}$, then the family of semi-norms $\{\sigma _{l}$%
, $l\in 
\mathbb{N}
\}$ such that%
\begin{equation}
\sigma _{l}(f)=\int_{K}\left\vert D^{l}f(y)\right\vert ^{2}dy)^{\frac{1}{2}},%
\text{ \ \ \ \ \ \ \ \ }f\in C^{\infty }(K)
\end{equation}%
define on $C^{\infty }(K)$ the same topology of the Frechet topology defined
by the semi-normas $\left\Vert X^{\alpha }f\right\Vert _{2}$ defined as%
\begin{equation}
\left\Vert X^{\alpha }f\right\Vert _{2}=\int_{K}(\left\vert X^{\alpha
}f(y)\right\vert ^{2}dy)^{\frac{1}{2}},\text{ \ \ \ \ \ \ \ \ }f\in
C^{\infty }(K)
\end{equation}%
where $\alpha =(\alpha _{1},$.....,$\alpha _{m})\in 
\mathbb{N}
^{m},$ see $[4,p.565]$

Let $\widehat{K}$ be the set of all irreducible unitary representations of $%
K.$ If $\gamma \in \widehat{K}$, we denote by $E_{\gamma }$ the space of the
representation $\gamma $ and $d_{\gamma }$\ its dimension then we get\qquad

\textbf{Definition 4.1.} \textit{The Fourier transform of a function }$f\in
C^{\infty }(K)$\textit{\ is defined as} 
\begin{equation}
Tf(\gamma )=\int\limits_{K}f(x)\gamma (x^{-1})dx
\end{equation}%
\textit{where }$T$\textit{\ is the Fourier transform on} $K$

\textbf{Theorem (\textbf{A. Cerezo }}$[4]$\textbf{) 4.1.} \textit{Let} $f\in
C^{\infty }(K),$ \textit{then we have the inversion of the Fourier transform}
\begin{equation}
f(x)=\sum\limits_{\gamma \in \widehat{K}}d\gamma tr[Tf(\gamma )\gamma (x)]
\end{equation}

\begin{equation}
f(I_{K})=\sum\limits_{\gamma \in \widehat{K}}d\gamma tr[Tf(\gamma )]
\end{equation}%
\textit{and the Plancherel formula} 
\begin{equation}
\left\Vert f(x)\right\Vert _{2}^{2}=\int_{K}\left\vert f(x)\right\vert
^{2}dx=\sum\limits_{\gamma \in \widehat{K}}d_{\gamma }\left\Vert Tf(\gamma
)\right\Vert _{H.S}^{2}
\end{equation}%
\textit{for any }$f\in L^{1}(K),$ \textit{where }$I_{K}$ \textit{is the
identity element of \ }$K$ \textit{and} $\left\Vert Tf(\gamma )\right\Vert
_{H.S}^{2}$ \textit{is the Hilbert- Schmidt norm of the operator} $Tf(\gamma
)$\ 

\bigskip \textbf{Definition 4.2}. \textit{For any function} $f\in \mathcal{D}%
(G),$ \textit{we can define a function} $\Upsilon (f)$\textit{on }$G\times K$
$=G\times SO(4)$ \textit{by} 
\begin{equation}
\Upsilon (f)(g,k_{1})=\Upsilon (f)(kna,k_{1})=f(gk_{1})=f(knak_{1})
\end{equation}%
\textit{for }$g=kna\in G,$ \textit{and} $k_{1}\in K$ . \textit{The
restriction of} $\ \Upsilon (f)\ast \psi (g,k_{1})$ \textit{on} $K(G)$ 
\textit{is }$\Upsilon (f)\ast \psi (g,k_{1})\downarrow
_{K(G)}=f(nak_{1})=f(g)\in \mathcal{D}(G),$ \textit{and }$\Upsilon
(f)(g,k_{1})\downarrow _{SO(4)}=f(kna)$ $\in \mathcal{D}(G)$

\textbf{Remark 4.1}. $\Upsilon (f)$ is invariant in the following sense%
\begin{equation}
\Upsilon (f)(gh,h^{-1}k_{1})=\Upsilon (f)(g,k_{1})
\end{equation}

\textbf{Definition 4.3}.\textit{\ If }$f$ \textit{and }$\psi $ \textit{are
two functions belong to} $\mathcal{D}(G),$ \textit{then we can define the
convolution of } $\Upsilon (f)\ $\textit{and} $\psi $\ \textit{on} $G\times
SO(4)$ \textit{as}

\begin{eqnarray*}
&&\Upsilon (f)\ast \psi (g,k_{1}) \\
&=&\int\limits_{G}\Upsilon (f)(gg_{2}^{-1},k_{1})\psi (g_{2})dg_{2} \\
&=&\int\limits_{SO(4)}\int\limits_{N}\int\limits_{A}\Upsilon
(f)(knaa_{2}^{-1}n_{2}^{-1}k^{-1}k_{1})\psi
(k_{2}n_{2}a_{2})dk_{2}dn_{2}da_{2}
\end{eqnarray*}%
and so we get 
\begin{eqnarray*}
\Upsilon (f)\ast \psi (g,k_{1}) &\downarrow &_{K(G)}=\Upsilon (f)\ast \psi
(I_{K}na,k_{1}) \\
&=&\int\limits_{SO(4)}\int\limits_{N}\int%
\limits_{A}f(naa_{2}^{-1}n_{2}^{-1}k^{-1}k_{1})\psi
(k_{2}n_{2}a_{2})dk_{2}dn_{2}da_{2} \\
&=&\Upsilon (f)\ast \psi (na,k_{1})
\end{eqnarray*}%
where $g_{2}=k_{2}n_{2}a_{2}$

\bigskip \textbf{Definition 4.3}. \textit{For }$f\in \mathcal{D}(G)$, 
\textit{let} $\Upsilon (f)$ \textit{be its associated function}, \textit{we
define the Fourier transform of \ }$\Upsilon (f)(g,k_{1})$ \textit{by }%
\begin{eqnarray}
&&\mathcal{F}\Upsilon (f))(I_{SO(4)},\xi ,\lambda ,\gamma )  \notag \\
&=&\int_{N}\int_{A}[\int_{SO(4)}(T\Upsilon (f)(I_{SO(4))}na,k_{1})\gamma
(k_{1}^{-1})dk_{1}]  \notag \\
&&a^{-i\lambda }e^{-\text{ }i\langle \text{ }\xi ,\text{ }n\rangle }\text{ }%
dadn  \notag \\
&=&\int\limits_{SO(3)}\int_{N}\int_{A}[\Upsilon
(f)(I_{SO(3)}na,k_{1})]a^{-i\lambda }e^{-\text{ }i\langle \text{ }\xi ,\text{
}n\rangle }\text{ }\gamma (k_{1}^{-1})dadndk_{1}
\end{eqnarray}%
\textit{where }$\mathcal{F}$ \textit{is the Fourier transform on }$AN$ 
\textit{and }$T$ \textit{is the Fourier transform on} $SO(4),$ $I_{SO(3)}$ 
\textit{is the identity element of }$SO(4)$, \textit{and }$%
n=(x_{1},x_{2},x_{3},x_{4},x_{5},x_{6}),n_{2}=(y_{1},y_{2},y_{3},y_{4},y_{5},y_{6}),\xi =(\xi _{1},\xi _{2},\xi _{3},\xi _{4},\xi _{5},\xi _{6}),a=b_{1}b_{2}b_{3} 
$ \textit{and} $a=a_{1}a_{2}a_{3}$

\textbf{Plancherel's Theorem 4.2\textit{\textbf{.} }}\textit{For any
function\ }$f\in $\textit{\ }$L^{1}(G)\cap $\textit{\ }$L^{2}(G),$\textit{we
get }%
\begin{equation}
\int_{G}\left\vert f(g)\right\vert
^{2}dg=\int\limits_{A}\int\limits_{N}\int\limits_{SO(4)}\left\vert
f(kna)\right\vert ^{2}dadndk=\sum\limits_{\gamma \in \widehat{SO(4)}%
}d_{\gamma }\int\limits_{\mathbb{R}^{3}}\int\limits_{\mathbb{R}%
^{6}}\left\Vert T\mathcal{F}f(\alpha ,\xi ,\gamma )\right\Vert
_{2}^{2}d\alpha d\xi
\end{equation}%
\textit{\ }%
\begin{equation}
f(I_{A}I_{N}I_{S^{1}})=\int\limits_{N}\int\limits_{A}\sum\limits_{\gamma \in 
\widehat{K}}d_{\gamma }T\mathcal{F}f(\alpha ,\xi ,\gamma )]d\alpha d\xi
=\sum\limits_{\gamma \in \widehat{K}}d_{\gamma }\int\limits_{\mathbb{R}%
^{3}}\int\limits_{\mathbb{R}^{6}}T\mathcal{F}f(\alpha ,\xi ,\gamma )d\alpha
d\xi
\end{equation}%
\textit{where }$,$ $I_{A},I_{N},$ and $I_{K}$ \textit{are the identity
elements of} $A$, $N$ \textit{and }$K$ \textit{respectively, }$\mathcal{F}$ 
\textit{is the Fourier transform on }$AN$ \textit{and }$T$ \textit{is the
Fourier transform on} $K,$

\bigskip \textit{Proof: }First let $\overset{\vee }{f}$ \ be the function
defined by 
\begin{equation}
\ \overset{\vee }{f}(kna)=\overline{f((kna)^{-1})}=\overline{%
f(a^{-1}n^{-1}k^{-1})}
\end{equation}

Then we have%
\begin{eqnarray}
&&\int_{SL(4,\mathbb{R})}\left\vert f(g)\right\vert ^{2}dg  \notag \\
&=&\Upsilon (f)\ast \overset{\vee }{f}(I_{SO(4)}I_{N}I_{A},I_{SO(4)})  \notag
\\
&=&\int\limits_{G}\Upsilon (f)(I_{SO(4)}I_{N}I_{A}g_{2}^{-1},I_{SO(4)})%
\overset{\vee }{f}(g_{2})dg_{2}  \notag \\
&=&\int\limits_{A}\int\limits_{N}\int_{SO(4)}\Upsilon
(f)(a_{2}^{-1}n_{2}^{-1}k_{2}^{-1},I_{SO(4)})\overset{\vee }{f}%
(k_{2}n_{2}a_{2})da_{2}dn_{2}dk_{2}  \notag \\
&=&\int\limits_{A}\int\limits_{N}%
\int_{SO(4)}f(a_{2}^{-1}n_{2}^{-1}k_{2}^{-1})\overline{%
f((k_{2}n_{2}a_{2})^{-1})}da_{2}dn_{2}dk_{2}  \notag \\
&=&\int\limits_{A}\int\limits_{N}\int_{SO(4)}\left\vert
f(a_{2}n_{2}k_{2})\right\vert ^{2}da_{2}dn_{2}dk_{2}
\end{eqnarray}

\bigskip Secondly%
\begin{eqnarray*}
&&\Upsilon (f)\ast \overset{\vee }{f}(I_{SO(4)}I_{N}I_{A},I_{SO(4)}) \\
&=&\int\limits_{\mathbb{R}^{18}}\sum\limits_{\gamma \in \widehat{SO(4)}%
}d\gamma \int_{SO(4)}tr(\Upsilon (f)\ast \overset{\vee }{f}%
(I_{SO(4)}na,k_{1})\gamma (k_{1}^{-1})) \\
&&a^{-i\alpha }e^{-\text{ }i\langle \text{ }\xi ,\text{ }n\rangle
}dadndk_{1}d\lambda d\xi \\
&=&\sum\limits_{\gamma \in \widehat{SO(4)}}d\gamma \int_{SO(4)}\int\limits_{%
\mathbb{R}^{18}}tr[\Upsilon (f)\ast \overset{\vee }{f}%
(I_{SO(3}na,k_{1})dka^{-i\alpha }e^{-\text{ }i\langle \text{ }\xi ,\text{ }%
n\rangle }\text{ }\gamma (k_{1}^{-1})]dadndk_{1}d\lambda d\xi \\
&=&\int\limits_{\mathbb{R}^{27}}\sum\limits_{\gamma \in \widehat{SO(4)}%
}\int_{SO(4)}tr[\Upsilon (f)(I_{SO(4)}nab^{-1}n_{2}^{-1}k_{2}^{-1},k_{1})%
\overset{\vee }{f}(k_{2}n_{2}b)\gamma (k_{1}^{-1})dk_{1}] \\
&&a^{-i\alpha }e^{-\text{ }i\langle \text{ }\xi ,\text{ }n\rangle
}dndadn_{2}dbd\lambda d\xi
\end{eqnarray*}%
where

\begin{eqnarray*}
&&e^{-i(y_{6}-x_{6}+x_{1}x_{5}-x_{1}y_{5}-x_{1}x_{2}x_{4}+x_{3}x_{4}-x_{3}y_{4}+x_{1}x_{2}y_{4})\xi _{6}}%
\text{ }e^{-i((y_{5}-x_{5}+x_{2}y_{4}-x_{2}x_{4})\xi _{5}+(y_{4}-x_{4})\xi
_{4})} \\
&&e^{-i((y_{3}-x_{3}-x_{1}y_{2}+x_{1}x_{2})\xi _{3}+(y_{2}-x_{2})\xi
_{2}+(y_{1}-x_{1})\xi _{1})}\text{ } \\
&=&e^{-i(y_{6}\xi _{6}-x_{6}\xi _{6}+y_{5}\xi _{5}-x_{5}\xi _{5}+y_{4}\xi
_{4}-x_{4}\xi _{4}+y_{3}\xi _{3}-\xi _{3}x_{3}-\xi _{2}x_{2}+y_{2}\xi
_{2}-\xi _{1}x_{1}+y_{1}\xi _{1})}
\end{eqnarray*}%
$%
n=(x_{1},x_{2},x_{3},x_{4},x_{5},x_{6}),n_{2}=(y_{1},y_{2},y_{3},y_{4},y_{5},y_{6}),\xi =(\xi _{1},\xi _{2},\xi _{3},\xi _{4},\xi _{5},\xi _{6}),a=b_{1}b_{2}b_{3} 
$ and $a=a_{1}a_{2}a_{3}$

\bigskip Using the fact that%
\begin{equation}
\int\limits_{A}\int\limits_{N}\int_{SO(4)}f(kna)dadndk=\int\limits_{N}\int%
\limits_{A}\int_{SO(4)}f(kan)a^{2}dndadk
\end{equation}%
and 
\begin{eqnarray}
&&\int\limits_{\mathbb{R}^{6}}\int\limits_{A}\int\limits_{N}%
\int_{SO(4)}f(kna)e^{-\text{ }i\langle \text{ }\xi ,\text{ }n\rangle
}dadndkd\xi  \notag \\
&=&\int\limits_{\mathbb{R}^{6}}\int\limits_{A}\int\limits_{N}%
\int_{SO(4)}f(kan)e^{-\text{ }i\langle \text{ }\xi ,\text{ }an_{1}a^{-1}%
\text{ }\rangle }a^{2}dadndkd\xi  \notag \\
&=&\int\limits_{\mathbb{R}^{6}}\int\limits_{A}\int\limits_{N}%
\int_{SO(4)}f(kan)e^{-\text{ }i\langle \text{ }a\xi a^{-1},\text{ }n\rangle
}a^{2}dadndkd\xi  \notag \\
&=&\int\limits_{\mathbb{R}^{6}}\int\limits_{A}\int\limits_{N}%
\int_{SO(4)}f(kan)e^{-\text{ }i\langle \text{ }\xi ,\text{ }n\rangle
}dadndkd\xi
\end{eqnarray}

Then we have 
\begin{eqnarray*}
&&\Upsilon (f)\ast \overset{\vee }{f}(I_{SO(4)}I_{N}I_{A},I_{SO(4)}) \\
&=&\int\limits_{\mathbb{R}27}\int_{SO(4)}\int\limits_{A}\int\limits_{N}\sum%
\limits_{\gamma \in \widehat{SO(3)}}d_{\gamma
}\int_{SO(4)}f(nab^{-1}n_{2}^{-1}k_{2}^{-1},k_{1})\overset{\vee }{f}%
(k_{2}n_{2}b)\gamma (k_{1}^{-1})dk_{1}dk_{2} \\
&&a^{-i\lambda }e^{-\text{ }i\langle \text{ }\xi ,\text{ }n\rangle
}dndadn_{2}da_{2}d\lambda d\xi \\
&=&\int\limits_{\mathbb{R}27}\sum\limits_{\gamma \in \widehat{SO(4)}%
}d_{\gamma }\int_{SO(4)}\int_{SO(4)}f(ab^{-1}nn_{2}^{-1}k_{2}^{-1},k_{1})%
\overset{\vee }{f}(k_{2}n_{2}b)\gamma (k_{1}^{-1})dk_{1}dk_{2} \\
&&a^{-i\lambda }e^{-\text{ }i\langle \text{ }\xi ,\text{ }n\rangle
}dndadn_{2}da_{2}d\lambda d\xi \\
&=&\int\limits_{\mathbb{R}27}\sum\limits_{\gamma \in \widehat{SO(4)}%
}d_{\gamma }\int_{SO(4)}\int_{SO(4)}f(ank_{2}^{-1},k_{1})\overset{\vee }{f}%
(k_{2}n_{2}b)\gamma (k_{1}^{-1})dk_{1}dk_{2} \\
&&ab^{-i\lambda }e^{-\text{ }i\langle \text{ }\xi ,\text{ }nn_{2}\rangle
}dndadn_{2}da_{2}d\lambda d\xi \\
&=&\int\limits_{\mathbb{R}27}\sum\limits_{\gamma \in \widehat{SO(4)}%
}d_{\gamma }\int_{SO(4)}\int_{SO(4)}f(ank_{2}^{-1}k_{1})\overset{\vee }{f}%
(k_{2}n_{2}b)\gamma (k_{1}^{-1})dk_{1}dk_{2} \\
&&ab^{-i\lambda }e^{-\text{ }i\langle \text{ }\xi ,\text{ }nn_{2}\rangle
}dndadn_{2}da_{2}d\lambda d\xi \\
&=&\int\limits_{\mathbb{R}27}\sum\limits_{\gamma \in \widehat{SO(4)}%
}d_{\gamma }\int_{SO(4)}\int_{SO(4)}f(ank_{1}^{-1})\overset{\vee }{f}%
(k_{2}n_{2}b)\gamma (k_{1}^{-1})\gamma (k_{2}^{-1})dk_{1}dk_{2} \\
&&a^{-i\lambda }b^{-i\lambda }e^{-\text{ }i\langle \text{ }\xi ,\text{ }%
n\rangle }e^{-\text{ }i\langle \text{ }\xi ,\text{ }n_{2}\rangle
}dndadn_{2}da_{2}d\lambda d\xi
\end{eqnarray*}

\bigskip So, we get%
\begin{eqnarray*}
&&\Upsilon (f)\ast \overset{\vee }{f}(I_{SO(4)}I_{N}I_{A},I_{SO(4)}) \\
&=&\int\limits_{\mathbb{R}27}\sum\limits_{\gamma \in \widehat{SO(4)}%
}d_{\gamma }\int_{SO(4)}\int_{SO(4)}f(ank_{1}^{-1})\overset{\vee }{f}%
(k_{2}n_{2}b)\gamma (k_{1}^{-1})\gamma (k_{2}^{-1})dk_{1}dk_{2} \\
&&a^{-i\lambda }b^{-i\lambda }e^{-\text{ }i\langle \text{ }\xi ,\text{ }%
n+n_{2}\rangle }dndadn_{2}da_{2}d\lambda d\xi \\
&=&\int\limits_{\mathbb{R}27}\sum\limits_{\gamma \in \widehat{SO(4)}%
}d_{\gamma }\int_{SO(4)}\int_{SO(4)}f(ank_{1}^{-1})\overline{%
f(b{}^{-1}n_{2}^{-1}k_{2}^{-1})}\gamma (k_{1}^{-1})\gamma
(k_{2}^{-1})dk_{1}dk_{2} \\
&&a^{-i\lambda }e^{-\text{ }i\langle \text{ }\xi ,\text{ }n\rangle
}b^{-i\lambda }e^{-\text{ }i\langle \text{ }\xi ,\text{ }n_{2}\rangle
}dndadn_{2}da_{2}d\lambda d\xi \\
&=&\int\limits_{\mathbb{R}27}\sum\limits_{\gamma \in \widehat{SO(4)}%
}d_{\gamma }\int_{SO(4)}\int_{SO(4)}f(ank_{1}^{-1})\overline{%
f(b{}n_{2}k_{2})\gamma (k_{2}^{-1})}\gamma (k_{1}^{-1})dk_{1}dk_{2} \\
&&a^{-i\lambda }e^{-\text{ }i\langle \text{ }\xi ,\text{ }n\rangle
}b^{-i\lambda }e^{\text{ }i\langle \text{ }\xi ,\text{ }n_{2}\rangle
}dndadn_{2}da_{2}d\lambda d\xi \\
&=&\int\limits_{\mathbb{R}27}\sum\limits_{\gamma \in \widehat{SO(4)}%
}d_{\gamma }\int_{SO(4)}\int_{SO(4)}f(ank_{1}^{-1})\overline{%
f(b{}n_{2}k_{2})\gamma (k_{2}^{-1})}\gamma (k_{1}^{-1})dk_{1}dk_{2} \\
&&a^{-i\lambda }e^{-\text{ }i\langle \text{ }\xi ,\text{ }n\rangle }%
\overline{b^{-i\lambda }e^{\text{ }-i\langle \text{ }\xi ,\text{ }%
n_{2}\rangle }}dndadn_{2}dbd\lambda d\xi \\
&=&\int\limits_{\mathbb{R}^{9}}\sum\limits_{\gamma \in \widehat{SO(4)}%
}d_{\gamma }T\mathcal{F}f(\lambda ,\xi ,\gamma )\overline{T\mathcal{F}%
f(\lambda ,\xi ,\gamma )}d\lambda d\xi =\int\limits_{\mathbb{R}%
^{9}}\sum\limits_{\gamma \in \widehat{SO(4)}}d_{\gamma }\left\vert T\mathcal{%
F}(f)(\lambda ,\xi ,\gamma )\right\vert ^{2}d\lambda d\xi
\end{eqnarray*}

Hence theorem of Plancherel on $SL(4,\mathbb{R})$ is Proved

Let $SP(4,\mathbb{R})=KNA$ be the Iwasawa decomposition of the symplectic $%
SP(4,\mathbb{R}).$ My state result is

\bigskip \textbf{Corollary\textit{\ }4.1.\textit{\ }}\textit{For any
function }$f\in $\textit{\ }$L^{1}(SP(4,\mathbb{R}))\cap $\textit{\ }$%
L^{2}(SP(4,\mathbb{R})),$\textit{we get}%
\begin{equation}
\int_{SP(4,\mathbb{R})}\left\vert f(v,g)\right\vert
^{2}dvdg=\int\limits_{N}\int\limits_{A}\sum_{\gamma \in \widehat{K}%
}d_{\gamma }\left\Vert \mathcal{F}_{\mathbb{R}^{2}}T\mathcal{F}F(\xi
,\lambda ,\gamma )\right\Vert ^{2}d\eta d\lambda d\xi
\end{equation}

Which is the Plancherel theorem on the symplectic $SP(4,\mathbb{R})$

\section{Plancherel Theorem on Group $\mathbb{R}^{4}\rtimes SL(4,\mathbb{R})$%
}

Let $P$ $=\mathbb{R}^{4}\rtimes _{\rho }SL(4,\mathbb{R})$ be the $14-$%
dimensional affine group$.$ Let $(v,g)$ and $(v^{\prime },g^{\prime })$ be
two elements belong $P,$ then the multiplication of $(v,g)$ and $(v^{\prime
},g^{\prime })$ is given by

\begin{equation}
(v,g)(v^{\prime },g^{\prime })=(v+\text{\ }\rho (g)(v^{\prime }),\text{ }%
gg^{\prime })=(v+\text{\ }gv^{\prime },\text{ }gg^{\prime })\text{ \ }
\end{equation}%
for any $(v,v^{\prime })\in \mathbb{R}^{4}\times \mathbb{R}^{4}$ and $(g,$ $%
g^{\prime })\in SL(4,\mathbb{R})\times SL(4,\mathbb{R}),$ where $gv^{\prime
}=\rho (g)(v^{\prime }).$ To define the Fourier transform on $P$, we
introduce the following new group

\textbf{Definition 5.1}. \textit{Let }$Q=\mathbb{R}^{3}$ $\times SL(4,%
\mathbb{R})\times SL(4,\mathbb{R})$ \textit{be the group with law}: 
\begin{eqnarray}
X\cdot Y &=&(v,h,g)(v^{\prime },h^{\prime },g^{\prime })  \notag \\
&=&(v+\text{\ }gv^{\prime },hh^{\prime },\text{ }gg^{\prime })
\end{eqnarray}%
\textit{for all} $X=(v,h,g)$ $\in Q$ \textit{and} $Y=(v^{\prime },h^{\prime
},g^{\prime })\in Q.$ Denote by $A=\mathbb{R}^{4}$ $\times SL(4,\mathbb{R})$
the group of the direct product of $\mathbb{R}^{4}$ with the group $SL(4,%
\mathbb{R}).$ Then the group $A$ can be regarded as the subgroup $\mathbb{R}%
^{3}$ $\times SL(4,\mathbb{R})\times \{I_{SL(4,\mathbb{R})}\}$ of $Q\mathbb{%
\ }$and $P$ can be regarded as the subgroup $\mathbb{R}^{4}$ $\times
\{I_{SL(4,\mathbb{R})}\}\times SL(4,\mathbb{R})$ of $Q.$

\textbf{Definition 5.2}. \textit{For any function} $f\in \mathcal{D}(P),$ 
\textit{we can define a function} $\widetilde{f}$ \textit{on }$Q$ \textit{by}

\begin{equation}
\widetilde{f}(v,g,h)=f(gv,gh)
\end{equation}

\textbf{Remark 5.1. }\textit{The function} $\widetilde{f}$ \textit{\ is
invariant in the following sense}%
\begin{equation}
\widetilde{f}(q^{-1}v,g,q^{-1}h)=\widetilde{f}(v,gq^{-1},h)\text{\ }
\end{equation}

\textbf{Theorem 5.1. }\textit{For any function }$\psi \in \mathcal{D}(P)$ 
\textit{and }$\widetilde{f}\in C^{\infty }(Q)$ \textit{invariant in sense }$%
(32)$\textit{, we get}%
\begin{equation}
\psi \ast \widetilde{f}(v,h,g)=\widetilde{f}\ast _{c}\psi (v,h,g)
\end{equation}%
\textit{\ where} $\ast $ \textit{signifies the convolution product on} $P$ 
\textit{with respect the variable} $(v,g),$ \textit{and }$\ast _{c}$\textit{%
signifies the convolution product on} $A$ \textit{with respect the variable} 
$(v,h)$ \ \ \ \ \ \ \ \ \ \ \ \ \ \ \ \ \ \ \ \ \ \ \ \ \ \ \ \ \ \ \ \ \ \
\ \ \ \ \ \ \ \ \ \ \ 

\textit{Proof : }In fact for each $\psi \in \mathcal{D}(P)$ \textit{and }$%
\widetilde{f}\in C^{\infty }(Q),$ we have%
\begin{eqnarray}
&&\psi \ast \widetilde{f}(v,h,g)  \notag \\
&=&\int\limits_{\mathbb{R}^{4}}\int_{SL(4,\mathbb{R})}\widetilde{f}%
((v^{\prime },g^{\prime })^{-1}(v,h,g))\psi (v^{\prime },g^{\prime
})dv^{\prime }dg^{\prime }  \notag \\
&=&\int\limits_{\mathbb{R}^{4}}\int_{SL(4,\mathbb{R})}\widetilde{f}%
[(g^{\prime }{}^{-1}(-v^{\prime }),g^{\prime }{}^{-1})(v,h,g)]\psi
(v^{\prime },g^{\prime })dv^{\prime }dg^{\prime }  \notag \\
&=&\int\limits_{\mathbb{R}^{4}}\int_{SL(4,\mathbb{R})}\widetilde{f}%
[(g^{\prime }{}^{-1}(-v^{\prime }),g^{\prime }{}^{-1})(v,h,g)]\psi
(v^{\prime },g^{\prime })dv^{\prime }dg^{\prime } \\
&=&\int\limits_{\mathbb{R}^{4}}\int_{SL(4,\mathbb{R})}\widetilde{f}%
[(g^{\prime }{}^{-1}(v-v^{\prime }),h,g^{\prime }{}^{-1}g)]\psi (v^{\prime
},g^{\prime })dv^{\prime }dg^{\prime }  \notag \\
&=&\int\limits_{\mathbb{R}^{4}}\int_{SL(4,\mathbb{R})}\widetilde{f}%
[v-v^{\prime },hg^{\prime }{}^{-1},g]\psi (v^{\prime },g^{\prime
})dv^{\prime }dg^{\prime }  \notag \\
&=&\widetilde{f}\ast _{c}\psi (v,h,g)
\end{eqnarray}

\textbf{Corollary 5.1.\ }\textit{From theorem \textbf{5.1}, the equation \
turns as} 
\begin{eqnarray}
&&\psi \ast \widetilde{f}(v,h,I_{G})  \notag \\
&=&\psi \ast _{c}\widetilde{f}(v,h,I_{SL(3,\mathbb{R})})=\int\limits_{%
\mathbb{R}^{2}}\int_{SL(3,\mathbb{R})}\widetilde{f}[v-v^{\prime },hg^{\prime
}{}^{-1},g]\psi (v^{\prime },g^{\prime })dv^{\prime }dg^{\prime }  \notag \\
&=&\int\limits_{\mathbb{R}^{4}}\int_{SL(4,\mathbb{R})}f[hg^{\prime
}{}^{-1}(v-v^{\prime }),hg^{\prime }{}^{-1}]\psi (v^{\prime },g^{\prime
})dv^{\prime }dg^{\prime }=h(f)\ast _{c}\psi (v,h)
\end{eqnarray}%
where

\begin{equation}
h(f)(v,g)=f(gv,g)
\end{equation}

\textbf{Definition 5}.\textbf{3}. \textit{Let }$\Upsilon F$ \textit{be the
function on }$P\times SL(4,\mathbb{R})$ \textit{defined by}%
\begin{equation}
\Upsilon F(v,(g,k_{1}))=F(v,gk_{1})
\end{equation}

\textbf{Definition 5.4. \ }\textit{Let }$\psi \in \mathcal{D}(P)$ \textit{%
and }$F\in \mathcal{D}(P),$\textit{then we can define a convolution product
on the Affine group }$P$ \textit{as}%
\begin{eqnarray*}
&&\psi \ast _{c}\Upsilon F(v,(g,k_{1})) \\
&=&\int\limits_{\mathbb{R}^{4}}\int_{SL(4,\mathbb{R})}\Upsilon F(v-v^{\prime
},(gg^{\prime }{}^{-1},k_{1}))\psi (v^{\prime },g^{\prime })dv^{\prime
}dg^{\prime } \\
&=&\int\limits_{\mathbb{R}^{4}}\int_{K}\int_{N}\int_{A}F(v-v^{\prime
},kna(k^{\prime }n^{\prime }a^{\prime })^{-1}k_{1}))\psi (v^{\prime
},k^{\prime }n^{\prime }a^{\prime })dv^{\prime }dk^{\prime }dn^{\prime
}da^{\prime }
\end{eqnarray*}%
\textit{where }$g=kna$ \textit{and} $g^{\prime }=k^{\prime }n^{\prime
}a^{\prime }$

\textbf{Corollary 5.2. }\textit{For any function }$F$ \textit{belongs to} $%
\mathcal{D}(P)$ , \textit{we obtain}%
\begin{eqnarray*}
&&\psi \ast _{c}\Upsilon h(F)(v,(g,k_{1})) \\
&=&\int\limits_{\mathbb{R}^{4}}\int_{SL(4,\mathbb{R})}\Upsilon
h(F)(v-v^{\prime },(gg^{\prime }{}^{-1},k_{1})\psi (v^{\prime },g^{\prime
})dv^{\prime }dg^{\prime } \\
&=&\int\limits_{\mathbb{R}^{4}}\int_{SL(4,\mathbb{R})}\Upsilon
h(F)(v-v^{\prime },(gg^{\prime }{}^{-1},k_{1})\psi (v^{\prime },g^{\prime
})dv^{\prime }dg^{\prime } \\
&=&\int\limits_{\mathbb{R}^{4}}\int_{SL(4,\mathbb{R})}h(F)(v-v^{\prime
},gg^{\prime }{}^{-1}k_{1})\psi (v^{\prime },g^{\prime })dv^{\prime
}dg^{\prime } \\
&=&\int\limits_{\mathbb{R}^{4}}\int_{SL(4,\mathbb{R})}F(gg^{\prime
}{}^{-1}k_{1}(v-v^{\prime }),gg^{\prime }{}^{-1}k_{1})\psi (v^{\prime
},g^{\prime })dv^{\prime }dg^{\prime }
\end{eqnarray*}

\textbf{Corollary 5.3. }\textit{For any function }$F$ \textit{belongs to} $%
\mathcal{D}(P)$ , \textit{we obtain}\textbf{\ } 
\begin{equation}
F\ast \Upsilon h(\overset{\vee }{F})(0,(I_{SL(4,\mathbb{R}%
)},I_{K}))=\int\limits_{\mathbb{R}^{4}}\int_{SL(4,\mathbb{R}%
)}|f(v,g)|^{2}dgdv=\left\Vert f\right\Vert _{2}^{2}
\end{equation}

\textit{Proof: }If\ $F\in \mathcal{D}(P),$then we get

\begin{eqnarray*}
&&F\ast \Upsilon h(\overset{\vee }{F})(0,(I_{SL(4,\mathbb{R})},I_{K})) \\
&=&\int\limits_{\mathbb{R}^{4}}\int_{SL(4,\mathbb{R})}\Upsilon \hbar (%
\overset{\vee }{F})[(0-v),(I_{G}g^{-1},I_{S^{1}})]F(v,g)dgdv \\
&=&\int\limits_{\mathbb{R}^{4}}\int_{SL(4,\mathbb{R})}\mathit{\ }\hbar (%
\overset{\vee }{F})[(0-v),I_{G}g^{-1}I_{S^{1}}]F(v,g)dgdv \\
&=&\int\limits_{SO(3)}\int\limits_{\mathbb{R}^{3}}\mathit{\ }\hbar (\overset{%
\vee }{F})[(-v),g^{-1}]F(v,g)dgdv \\
&=&\int\limits_{\mathbb{R}^{4}}\int_{SL(4,\mathbb{R})}\overset{\vee }{F}%
[g^{-1}(-v),g^{-1}]F(v,g)dgdv \\
&=&\int\limits_{\mathbb{R}^{4}}\int_{SL(4,\mathbb{R})}\overline{%
F[g^{-1}(-v),g^{-1}]^{-1}}F(v,g)dgdv \\
&=&\int\limits_{\mathbb{R}^{4}}\int_{SL(4,\mathbb{R})}\overline{F[v,g]}%
F(v,g)dgdv=\int\limits_{\mathbb{R}^{4}}\int_{SL(4,\mathbb{R})}\left\vert
f(v,g)\right\vert ^{2}dgdv
\end{eqnarray*}

\textbf{Definition 5.5.}\textit{\ Let} $f\in \mathcal{D}(P),$ \textit{we
define its Fourier transform by}%
\begin{equation}
\mathcal{F}_{\mathbb{R}^{4}}T\mathcal{F}f(\eta ,\gamma ,\xi ,\lambda
)=\int\limits_{\mathbb{R}^{4}}\int_{A}\int_{N}\int_{K}f(v,kna)e^{-\text{ }%
i\langle \text{ }\eta ,\text{ }v\rangle }\text{ }\gamma (k^{-1})a^{-i\lambda
}e^{-\text{ }i\langle \text{ }\xi ,\text{ }n\rangle }dkdadndv  \notag
\end{equation}%
\textit{where} $\mathcal{F}_{\mathbb{R}^{4}}$ \textit{is the Fourier
transform on} $\mathbb{R}^{3},$ $kna=g,$ $\eta =(\eta _{1},\eta _{2},\eta
_{3},\eta _{4})\in \mathbb{R}^{4},$ $v=(v_{1},v_{2},v_{3},v_{4})\in \mathbb{R%
}^{4},$ \textit{and} $dv=dv_{1}dv_{2}dv_{3}dv_{4}$ \textit{is the Lebesgue
measure on} $\mathbb{R}^{4}$ and

\begin{equation}
\langle (\eta _{1},\eta _{2},\eta _{3},\eta
_{4}),(v_{1},v_{2},v_{3},v_{4})\rangle =\sum_{i=1}^{4}\eta _{i}v_{i}
\end{equation}

\textbf{Plancherel's Theorem 5.2\textit{\textbf{\textit{.}} }}\textit{For
any function }$f\in $\textit{\ }$L^{1}(P)\cap $\textit{\ }$L^{2}(P),$\textit{%
we get}%
\begin{equation}
\int_{P}\left\vert f(v,g)\right\vert ^{2}dvdg=\int\limits_{\mathbb{R}%
^{13}}\sum_{\gamma \in \widehat{SO(4)}}d_{\gamma }\left\Vert \mathcal{F}_{%
\mathbb{R}^{2}}T\mathcal{F}F(\eta ,\gamma ,\xi ,\lambda )\right\Vert
^{2}d\eta d\lambda d\xi
\end{equation}

\textit{Proof: }Let\ $\Upsilon \hbar (\overset{\vee }{F})$ be the function
defined as\textit{\ }%
\begin{eqnarray}
&&\mathit{\ }\Upsilon \hbar (\overset{\vee }{F})\text{\ }(v;(g,k_{1}))=%
\mathit{\ }\hbar (\overset{\vee }{F})\text{\ }(v;gk_{1})  \notag \\
&=&\overset{\vee }{F}(gk_{1}v;gk_{1})=\mathit{\ }\overline{%
F(gk_{1}v;gk_{1})^{-1})}
\end{eqnarray}%
then, we have

\begin{eqnarray*}
&&F\ast \Upsilon \mathit{\ }\hbar (\overset{\vee }{F})\text{\ }(0,(I_{SL(4,%
\mathbb{R})}I_{N}I_{A},I_{SO(4)})) \\
&=&\int\limits_{\mathbb{R}^{2}}\int\limits_{\mathbb{R}^{3}}\int\limits_{%
\mathbb{R}^{3}}\mathcal{F}_{\mathbb{R}^{4}}\mathcal{F(}F\ast \Upsilon 
\mathit{\ }\hbar (\overset{\vee }{F}\text{\ })(\eta ,(I_{SO(4)}\xi ,\lambda
,I_{SO(4)}))d\lambda d\xi d\eta \\
&=&\int\limits_{\mathbb{R}^{26}}\sum_{\gamma \in \widehat{SO(4)}}d_{\gamma
}\int\limits_{SO(4)}\mathcal{F}_{\mathbb{R}^{4}}T\mathcal{F(}F\ast \Upsilon 
\mathit{\ }\hbar (\overset{\vee }{F}))((\eta ,(I_{SO(4)}na,k_{1}))\gamma
(k_{1}^{-1})dk_{1}) \\
&&e^{-\text{ }i\langle \text{ }\eta ,\text{ }v\rangle }a^{-i\lambda }e^{-%
\text{ }i\langle \text{ }\xi ,\text{ }n\rangle }dadndvd\lambda d\xi d\eta \\
&=&\int\limits_{SL(4,\mathbb{R})}\int\limits_{\mathbb{R}^{30}}\sum_{\gamma
\in \widehat{SO(3)}}d_{\gamma }\int\limits_{SO(4)}\Upsilon \mathit{\ }\hbar (%
\overset{\vee }{F})((v-w),(I_{SO(4)}nag_{2}^{-1},k_{1}))\gamma
(k_{1}^{-1})dk_{1}F(w,g_{2})dg_{2} \\
&&e^{-\text{ }i\langle \text{ }\eta ,\text{ }v\rangle }a^{-i\lambda }e^{-%
\text{ }i\langle \text{ }\xi ,\text{ }n\rangle }dadndvdwd\lambda d\xi d\eta
\\
&=&\int\limits_{\mathbb{R}^{40}}\int\limits_{SO(4)}\sum_{\gamma \in \widehat{%
SO(4)}}d_{\gamma }\int\limits_{SO(4)}\hbar (\overset{\vee }{F}%
)((v-w),(nab^{-1}n_{2}^{-1}k_{2}^{-1}k_{1}))\gamma (k_{1}^{-1})dk_{1} \\
&&F(w,k_{2}n_{2}b)dk_{2}e^{-\text{ }i\langle \text{ }\eta ,\text{ }v\rangle
}a^{-i\lambda }e^{-\text{ }i\langle \text{ }\xi ,\text{ }n\rangle
}dbdn_{2}dadndwdvd\lambda d\xi d\eta \\
&=&\int\limits_{\mathbb{R}^{40}}\int\limits_{SO(4)}\sum_{\gamma \in \widehat{%
SO(4)}}d_{\gamma }\int\limits_{SO(4)}\hbar (\overset{\vee }{F}%
)(v,(ank_{1}))\gamma (k_{1}^{-1})dk_{1}F(w,k_{2}n_{2}b)\gamma
(k_{2}^{-1})dk_{2} \\
&&e^{-\text{ }i\langle \text{ }\eta ,\text{ }v+w\rangle }a^{-i\lambda
}b^{-i\lambda }e^{-\text{ }i\langle \text{ }\xi ,\text{ }n\rangle }e^{-\text{
}i\langle \text{ }\xi ,\text{ }n_{2}\rangle }da_{2}dn_{2}dadndwdvd\lambda
d\xi d\eta
\end{eqnarray*}

\bigskip So, we get%
\begin{eqnarray*}
&=&\int\limits_{\mathbb{R}^{40}}\int\limits_{SO(4)}\sum_{\gamma \in \widehat{%
SO(4)}}d_{\gamma }\int\limits_{SO(4)}(\overset{\vee }{F})(ank_{1}v,ank_{1})%
\gamma (k_{1}^{-1})dk_{1}F(w,k_{2}n_{2}a_{2})\gamma (k_{2}^{-1})dk_{2} \\
&&e^{-\text{ }i\langle \text{ }\eta ,\text{ }v\rangle }e^{-\text{ }i\langle 
\text{ }\eta ,\text{ }w\rangle }e^{-i\langle \lambda ,a\rangle }e^{-i\langle
\lambda ,a_{2}\rangle }e^{-\text{ }i\langle \text{ }\xi ,\text{ }n\rangle
}e^{-\text{ }i\langle \text{ }\xi ,\text{ }n_{2}\rangle
}da_{2}dn_{2}dadndwdvd\lambda d\xi d\eta \\
&=&\int\limits_{\mathbb{R}^{40}}\int\limits_{SO(4)}\sum_{\gamma \in \widehat{%
SO(4)}}d_{\gamma }\int\limits_{SO(4)}\overline{F(-v,k_{1}^{-1}n^{-1}a^{-1})}%
F(w,k_{2}n_{2}a_{2})\gamma (k_{1}^{-1})\gamma (k_{2}^{-1})\gamma (k_{2}^{-1})
\\
&&e^{-\text{ }i\langle \text{ }\eta ,\text{ }v\rangle }e^{-\text{ }i\langle 
\text{ }\eta ,\text{ }w\rangle }a^{-i\lambda }a_{2}^{-i\lambda }e^{-\text{ }%
i\langle \text{ }\xi ,\text{ }n\rangle }e^{-\text{ }i\langle \text{ }\xi ,%
\text{ }n_{2}\rangle }dk_{1}dk_{2}da_{2}dn_{2}dadndwdvd\lambda d\xi d\eta \\
&=&\int\limits_{\mathbb{R}0}\int\limits_{SO(4)}\sum_{\gamma \in \widehat{%
SO(4)}}d_{\gamma }\int\limits_{SO(4)}\overline{F(v,k_{1}na)\gamma (k_{1})}%
F(w,k_{2}n_{2}a_{2})\gamma (k_{2}^{-1})dk_{1}dk_{2} \\
&&e^{\text{ }i\langle \text{ }\eta ,\text{ }v\rangle }e^{-\text{ }i\langle 
\text{ }\eta ,\text{ }w\rangle }e^{i\langle \lambda ,a\rangle }e^{-i\langle
\lambda ,a_{2}\rangle }e^{\text{ }i\langle \text{ }\xi ,\text{ }n\rangle
}e^{-\text{ }i\langle \text{ }\xi ,\text{ }n_{2}\rangle
}da_{2}dn_{2}dadndwdvd\lambda d\xi d\eta \\
&&\int\limits_{\mathbb{R}^{13}}\sum_{\gamma \in \widehat{SO(4)}}d_{\gamma
}\left\Vert \mathcal{F}_{\mathbb{R}^{4}}T\mathcal{F}F(\eta ,\gamma ,\xi
,\lambda )\right\Vert _{H.S}^{2}d\eta d\lambda d\xi
\end{eqnarray*}

\bigskip Hence the theorem is proved on the $\mathbb{R}^{4}$ $\rtimes SL(4,%
\mathbb{R}).$

\textbf{Corollary 5.3.\textit{\ }}\textit{For any function }$f\in $\textit{\ 
}$L^{1}(\mathbb{R}^{4}\rtimes _{\rho }SP(4,\mathbb{R}))\cap $\textit{\ }$%
L^{2}(\mathbb{R}^{4}\rtimes _{\rho }SP(4,\mathbb{R})),$\textit{we get}%
\begin{equation}
\int_{\mathbb{R}^{4}\rtimes _{\rho }SP(4,\mathbb{R})}\left\vert
f(v,g)\right\vert ^{2}dvdg=\int\limits_{\mathbb{R}^{4}}\int_{N}\int_{A}%
\sum_{\gamma \in \widehat{K}}d_{\gamma }\left\Vert \mathcal{F}_{\mathbb{R}%
^{2}}T\mathcal{F}F(\eta ,\gamma ,\xi ,\lambda )\right\Vert ^{2}d\eta
d\lambda d\xi
\end{equation}

Which is the Plancherel theorem on the inhomogeneous group $\mathbb{R}%
^{4}\rtimes _{\rho }SP(4,\mathbb{R})$ of the symplectic $SP(4,\mathbb{R}),$
where $KNA$ is the Iwasawa decomposition of the symplectic group $SP(4,%
\mathbb{R}).$

\section{Hypoellipticity of Differential Operators on the Symplectic}

\textbf{6.1.} Denote by $SP_{N}$ the nilpotent symplectic subgroup of the
group $SP(4,\mathbb{R)}$ consists of all matrices of the form 
\begin{equation}
SP_{N}=\left\{ \left( 
\begin{array}{cccc}
1 & x & y & z \\ 
0 & 1 & z-xt & t \\ 
0 & 0 & 1 & 0 \\ 
0 & 0 & -x & 1%
\end{array}%
\right) ,(x,y,z,t)\in \mathbb{R}^{4}\right\}
\end{equation}

\bigskip

We denote by $N$ the nilpotent symplectic subgroup of $SP_{N}$, formed by
the following matrix

\begin{equation}
N=\left\{ \left( 
\begin{array}{cccc}
1 & x & y & z \\ 
0 & 1 & z & 0 \\ 
0 & 0 & 1 & 0 \\ 
0 & 0 & -x & 1%
\end{array}%
\right) ,(x,y,z)\in \mathbb{R}^{3}\right\}
\end{equation}

\ \ \ The group $N$ is isomorphic onto the group $G=\mathbb{R}^{3}\rtimes
_{\rho }\mathbb{R}$ semidirect of two groups $\mathbb{R}^{2}$ and $\mathbb{R}
$ $,$ where $\rho :\mathbb{R\rightarrow }$ $Aut(\mathbb{R}^{2})$ is the
group homomorphism defined by \ $\rho (x)(z,y)=(z+xy,y).$\ The
multiplication of two elements $X=(z,y,x)$ and $Y=(c,b,a)$ is given by

\begin{eqnarray}
&&(z,y,x)(c,b,a)  \notag \\
&=&(z+c+xb-ay,y+b,x+a)  \notag \\
&=&(z+c+\left\vert 
\begin{array}{c}
x\text{ \ \ }y \\ 
a\text{ \ \ }b%
\end{array}%
\right\vert ,y+b,x+a)
\end{eqnarray}%
$.$Our aim is to prove the solvability and hypoellipticity of the following
Lewy operators%
\begin{equation}
L=(-\partial _{x}-i\partial _{y}-2y\partial _{z}+2ix\partial _{z})
\end{equation}

\begin{equation}
L_{\star }=(-\partial _{x}+i\partial _{y}-2y\partial _{z}-2ix\partial _{z})
\end{equation}

\textbf{Definition 6.1. }\textit{One can define a transformation }$\hbar :$%
\textit{\ }$\mathcal{D}^{\prime }(\mathbb{R}^{3})\rightarrow \mathcal{D}%
^{\prime }(\mathbb{R}^{3})$

\begin{equation}
\hbar \Psi (z,y\text{ },x)=\Psi (z-2xy\text{ },y\text{ ,}-x)
\end{equation}

It results from this definition that $\hbar ^{2}=\hbar $

\bigskip \textbf{Theorem 6.1. }\textit{Let }$Q=\partial _{x}-i\partial _{y}$ 
\textit{be the Cauchy}-\textit{Riemann operator, then we have for any }$f\in
C^{\infty }(\mathbb{R}^{3})$%
\begin{equation}
(Lf)(z,y,-x)=\hbar Q\hbar f(z,y,-x)\ 
\end{equation}%
\newline

For the proof of this theorem see $[6].$

\textbf{Corollary 6.1. }\textit{The Lewy operator }$L$\textit{\ is solvable}

\textit{Proof: }In fact the Cauchy-Riemann operator $Q=\partial
_{x}-i\partial _{y}$ is solvable, because $QC^{\infty }(\mathbb{R}%
^{3})=C^{\infty }(\mathbb{R}^{3}),$ and $\hbar C^{\infty }(\mathbb{R}%
^{3})=C^{\infty }(\mathbb{R}^{3}).$ So, I have $LC^{\infty }(H)=C^{\infty
}(H)$

\textbf{Definition 6.1. }\textit{Let }$G$ \textit{be a Lie group an}\textbf{%
\ }\textit{operator \ }$\Gamma :\mathcal{D}^{\prime }(G)\rightarrow \mathcal{%
D}^{\prime }(G)$ \textit{is called hypoelliptic if }%
\begin{equation}
\Gamma \varphi \in C^{\infty }(G)\Longrightarrow \varphi \in C^{\infty }(G)
\end{equation}%
\textit{\ for every distribution} $\varphi \in \mathcal{D}^{\prime }(G).$

\bigskip \textbf{Theorem 6.2. }\textit{The Lewy operator is hypoelliptic }

\textit{Proof: }First\textit{\ }the operator\textit{\ }$\hbar $ is
hypoelliptic, and the Cauchy- Riemann operator $\partial _{x}-i\partial _{y}$
is hypoelliptic. So if $\varphi \in \mathcal{D}^{\prime }(\mathbb{R}^{3})%
\mathcal{\ }$and if $L\varphi (z,y,-x)=\hbar Q\hbar \varphi (z,y,-x)\in
C^{\infty }(\mathbb{R}^{3}),$ then I get 
\begin{eqnarray}
L\varphi &\in &C^{\infty }(\mathbb{R}^{3})\Longrightarrow \hbar Q\hbar
\varphi \in C^{\infty }(\mathbb{R}^{3})  \notag \\
&\Longrightarrow &Q\hbar \varphi \in C^{\infty }(\mathbb{R}^{3})\text{ }%
\Longrightarrow \hbar \varphi \in C^{\infty }(\mathbb{R}^{3})  \notag \\
&\Longrightarrow &\varphi \in C^{\infty }(\mathbb{R}^{3})
\end{eqnarray}

\textbf{Theorem 6.3.}\textit{\ Let }$Q_{\star }$\textit{\ be the operator}%
\begin{equation}
L_{\star }=(-\partial _{x}+i\partial _{y}-2y\partial _{z}-2ix\partial _{z})
\end{equation}%
\begin{equation}
Q_{\star }=\partial _{x}+i\partial _{y}
\end{equation}%
\textit{then for every }$\varphi \in C^{\infty }(\mathbb{R}^{3}),$\textit{\
I have}%
\begin{eqnarray}
&&\hbar (\partial _{x}-i\partial _{y})(\partial _{x}+i\partial _{y})\hbar
\varphi (z,y,-x)=\hbar \Delta \hbar \varphi (z,y,-x)  \notag \\
&=&[(-\partial _{x}-2y\partial _{z})+(-i\partial _{y}+2ix\partial
_{z})((-\partial _{x}-2y\partial _{z})+(i\partial _{y}-2ix\partial
_{z})\varphi ](z,y,-x)  \notag \\
&=&LL_{\star }\varphi (z,y,-x)
\end{eqnarray}%
where $\Delta $ and $L_{\star }$ are the operators%
\begin{equation}
\Delta =\frac{\partial ^{2}}{\partial _{x^{2}}}+\frac{\partial ^{2}}{%
\partial _{y^{2}}}
\end{equation}%
\begin{equation}
L_{\star }=(i\partial _{y}-2ix\partial _{z})+(-\partial _{x}-2y\partial _{z})
\end{equation}%
$L_{\star }$ is called the conjugate of the Lewy operator, which can be
considered another form of the Lewy operator. As in theorem \textbf{6.2,} we
can easily see that $L_{\star }C^{\infty }(\mathbb{R}^{3})=C^{\infty }(%
\mathbb{R}^{3}).$ The operator $LL_{\star }$ can be regarded as the square
of the Lewy operator on the $3-$dimensional Heisenberg group.

\textbf{Corollary 6.1. }\textit{The operators }$LL_{\star }$ \textit{and }$%
L_{\star }$ \textit{are hypoelliptic}

\textit{Proof: }From the above we deduce the following%
\begin{eqnarray}
L_{\star }\varphi &\in &C^{\infty }(\mathbb{R}^{3})\Longrightarrow \hbar
Q_{\star }\hbar \varphi \in C^{\infty }(\mathbb{R}^{3})  \notag \\
&\Longrightarrow &Q_{\star }\hbar \varphi \in C^{\infty }(\mathbb{R}^{3})%
\text{ }\Longrightarrow \hbar \varphi \in C^{\infty }(\mathbb{R}^{3})  \notag
\\
&\Longrightarrow &\varphi \in C^{\infty }(\mathbb{R}^{3})
\end{eqnarray}

In other hand we have%
\begin{eqnarray}
LL_{\star }\varphi &\in &C^{\infty }(\mathbb{R}^{3})\Longrightarrow \hbar
QQ_{\star }\hbar \varphi \in C^{\infty }(\mathbb{R}^{3})  \notag \\
&\Longrightarrow &QQ_{\star }\hbar \varphi \in C^{\infty }(\mathbb{R}^{3})%
\text{ }\Longrightarrow Q_{\star }\hbar \varphi \in C^{\infty }(\mathbb{R}%
^{3})  \notag \\
&\Longrightarrow &\hbar \varphi \in C^{\infty }(\mathbb{R}%
^{3})\Longrightarrow \varphi \in C^{\infty }(\mathbb{R}^{3})
\end{eqnarray}

\bigskip \textbf{Theorem 6.4.}\textit{\ The following left invariant
differential operators on }$G$%
\begin{equation}
y\partial _{z}+\partial _{x}+i\partial _{y}+ix\partial _{z}
\end{equation}%
\begin{equation}
\frac{\partial ^{2}}{\partial _{x^{2}}}-\frac{\partial ^{2}}{\partial
_{y^{2}}}-2x\frac{\partial }{\partial _{z}}\frac{\partial }{\partial y}+2y%
\frac{\partial }{\partial _{z}}\frac{\partial }{\partial x}+(y^{2}-x^{2})%
\frac{\partial ^{2}}{\partial _{z^{2}}}+\frac{\partial ^{2}}{\partial
_{z^{2}}}
\end{equation}%
\textit{are solvable and hypoelliptic}

\bigskip \textit{Proof: }The solvability results from theorem \textbf{6.1.}
For the hypoellipticity, we consider the mapping $\Gamma :\mathcal{D}%
^{\prime }(G)\rightarrow \mathcal{D}^{\prime }(G)$ defined by%
\begin{equation}
\Gamma \phi (z,y,x)=\phi (z-xy,y,x)
\end{equation}

The operator $\Gamma $ is hypoelliptic and its inverse is 
\begin{equation}
\Gamma ^{-1}\phi (z,y,x)=\phi (z+xy,y,x)
\end{equation}%
thus we get%
\begin{equation}
\Gamma (\partial _{x}+i\partial _{y})\Gamma ^{-1}\phi (z,y,x)=(y\partial
_{z}+\partial _{x}+i\partial _{y}+ix\partial _{z})\phi (z,y,x
\end{equation}%
and

\begin{eqnarray}
&&\Gamma (\frac{\partial ^{2}}{\partial _{x^{2}}}+\frac{\partial ^{2}}{%
\partial _{y^{2}}}+\frac{\partial ^{2}}{\partial _{z^{2}}})\Gamma ^{-1}\phi
(z,y,x) \\
&=&(\frac{\partial ^{2}}{\partial _{x^{2}}}-\frac{\partial ^{2}}{\partial
_{y^{2}}}-2x\frac{\partial }{\partial _{z}}\frac{\partial }{\partial y}+2y%
\frac{\partial }{\partial _{z}}\frac{\partial }{\partial x}+(y^{2}-x^{2})%
\frac{\partial ^{2}}{\partial _{z^{2}}}+\frac{\partial ^{2}}{\partial
_{z^{2}}})\phi (z,y,x)
\end{eqnarray}

Since the operators $\Gamma ,\partial _{x}+i\partial _{y},\frac{\partial ^{2}%
}{\partial _{x^{2}}}+\frac{\partial ^{2}}{\partial _{y^{2}}}+\frac{\partial
^{2}}{\partial _{z^{2}}}$ and $\Gamma ^{-1}$ are hypoelliptic, then the
hypoellipticity of the operators $(y\partial _{z}+\partial _{x}+i\partial
_{y}+ix\partial _{z})$ and $\frac{\partial ^{2}}{\partial _{x^{2}}}-\frac{%
\partial ^{2}}{\partial _{y^{2}}}-2x\frac{\partial }{\partial _{z}}\frac{%
\partial }{\partial y}+2y\frac{\partial }{\partial _{z}}\frac{\partial }{%
\partial x}+(y^{2}-x^{2})\frac{\partial ^{2}}{\partial _{z^{2}}}+\frac{%
\partial ^{2}}{\partial _{z^{2}}}$ is fulfilled

\bigskip \textbf{Hormander condition for the hypoellipticity }

By the sufficient condition of the hypoellipticity given by the Hormander
theorem $[3,$ page $11]$, we oblige already quoted the sublaplacian%
\begin{equation}
\frac{\partial ^{2}}{\partial _{x^{2}}}+\frac{\partial ^{2}}{\partial
_{y^{2}}}+4x\frac{\partial }{\partial _{z}}\frac{\partial }{\partial y}-4y%
\frac{\partial }{\partial _{z}}\frac{\partial }{\partial x}+4(y^{2}+x^{2})%
\frac{\partial ^{2}}{\partial _{z^{2}}}
\end{equation}%
which is hypoelliptic by the Hormander theorem, while the operator 
\begin{equation}
\frac{\partial ^{2}}{\partial _{x^{2}}}+\frac{\partial ^{2}}{\partial
_{y^{2}}}+4x\frac{\partial }{\partial _{z}}\frac{\partial }{\partial y}-4y%
\frac{\partial }{\partial _{z}}\frac{\partial }{\partial x}+4(y^{2}+x^{2})%
\frac{\partial ^{2}}{\partial _{z^{2}}}-4i\frac{\partial }{\partial _{z}}
\end{equation}%
is not hypoelliptic because the Hormander condition is not fulfilled. By
contrast all our results, which are obtained by above theorems, contradict
the Hormonder conditions for the solvability and the hypoellipticity.

The basis of the Lie algebra of the group $N$ is given by the following
vector fields $Z=\frac{\partial }{\partial _{z}},Y=(x\frac{\partial }{%
\partial _{z}}+\frac{\partial }{\partial y}),$ $X=(-y\frac{\partial }{%
\partial _{z}}+\frac{\partial }{\partial x})$. Since $[X,Y]=2Z,$ and $%
X,Y,[X,Y]$ span the Lie algebra of $N$ . Then the Hormander theorem in $[5]$%
, gives the hypoellipticity of the operator%
\begin{equation}
X^{2}+Y^{2}=(x\frac{\partial }{\partial _{z}}+\frac{\partial }{\partial y}%
)^{2}+(-y\frac{\partial }{\partial _{z}}+\frac{\partial }{\partial x})^{2}
\end{equation}

While my results prove the solvability and hypoellipticity operators

\begin{equation}
X^{2}+Y^{2}+Z^{2}=(x\frac{\partial }{\partial _{z}}+\frac{\partial }{%
\partial y})^{2}+(-y\frac{\partial }{\partial _{z}}+\frac{\partial }{%
\partial x})^{2}+\frac{\partial ^{2}}{\partial _{z^{2}}}
\end{equation}

As well known the Laplace operator 
\begin{equation}
\Delta =\dsum\limits_{i=1}^{3}\frac{\partial ^{2}}{\partial _{x_{i}^{2}}}
\end{equation}%
on the real vector group $\mathbb{R}^{3}$ is solvable and hypoelliptic. This
operator as a left invariant differential on the group $N$ is nothing but
the following operator%
\begin{equation}
\Delta _{h_{1}}=\dsum\limits_{i=1}^{3}\frac{\partial ^{2}}{\partial _{z^{2}}}%
+(x\frac{\partial }{\partial _{z}}+\frac{\partial }{\partial _{y}})^{2}+(-y%
\frac{\partial }{\partial _{z}}+\frac{\partial }{\partial _{x}})^{2}
\end{equation}%
and as a right invariant on $N$ is the operator 
\begin{equation}
\Delta _{h_{2}}=\dsum\limits_{i=1}^{3}\frac{\partial ^{2}}{\partial _{z^{2}}}%
+(-x\frac{\partial }{\partial _{z}}+\frac{\partial }{\partial _{y}})^{2}+(y%
\frac{\partial }{\partial _{z}}+\frac{\partial }{\partial _{x}})^{2}
\end{equation}%
where $\Delta _{h_{1}}(resp.\Delta _{h_{2}})$ is the left (resp. right)
invariant differential operator associated to $\Delta .$ The operators $%
\Delta _{h_{1}}$ and $\Delta _{h_{2}}$ can be regarded as the Laplacian
operators on the $3-$dimensional Symplectic Nilpotent group $N$.

My aim result is

\textbf{Theorem 6.4. }\textit{The} \textit{Laplace operators }$\Delta
_{h_{1}}$ \textit{and} $\Delta _{h_{2}}$ \textit{on the Heisenberg group are
hypoelliptic}

\textit{Proof: }We consider the following mappings from $\mathcal{D}^{\prime
}(N)\rightarrow \mathcal{D}^{\prime }(N)$ defined by%
\begin{equation}
\Lambda \Psi (z,y\text{ },x)=\Psi (z+xy\text{ },y\text{ ,}x)
\end{equation}

\begin{equation}
\tau \Psi (z,y\text{ },x)=\Psi (z+xy\text{ },-y\text{ ,}x)
\end{equation}%
\begin{equation}
\pi \Psi (z,y\text{ },x)=\Psi (z+xy\text{ },y\text{ ,}-x)
\end{equation}

These operators has the property of hypoellipticity, because if $\Lambda
\Psi (z,y,x)=\Psi (z+xy,-y,x)\in C^{\infty }(N),$ then $\Psi (z,y,x)\in
C^{\infty }(N)$, so on $\tau $ and $\pi .$ In other side we have%
\begin{equation}
\tau \Delta \Lambda \Psi (z,y,x)=\Delta _{h_{1}}\Psi (z,-y,x)
\end{equation}%
\begin{equation}
\pi \Delta \Lambda \Psi (z,y,x)=\Delta _{h_{2}}\Psi (z,y,-x)
\end{equation}

Since $\Delta $, $\Lambda ,\tau $ and $\pi $ are hypoelliptic, then the
hypoellipticity of $\Delta _{h_{1}}$and $\Delta _{h_{\substack{ 2}}}$are
accomplished

\section{\protect\bigskip On the Existence Theorem on $N$}

Out of the proofs of my book $[6]$, I solve here by different method the
equation%
\begin{equation*}
PC^{\infty }(G)=C^{\infty }(G)
\end{equation*}

For this, I introduce two groups: The first is the group $G\times \mathbb{R}%
, $ which is the direct product of the group $G$\ with the real vector group 
$\mathbb{R}.$ The second is the group $E=\mathbb{R}^{2}\times \mathbb{R}%
\times \mathbb{R}$ with law:

\begin{equation}
g\cdot g^{\prime }=(X,x,y)(X^{\prime },x^{\prime },y^{\prime })=(X+X^{\prime
}+yX^{\prime },,x+x^{\prime },y+y^{\prime })
\end{equation}

for all $g=(X,x,y)\in \mathbb{R}^{4},g^{\prime }=(X^{\prime },x^{\prime
},y^{\prime })\in \mathbb{R}^{4},X\in \mathbb{R}^{2}$ and $X^{\prime }\in 
\mathbb{R}^{2}$. In this case the group $G$ can be identified with the
closed sub$-$group $\mathbb{R}^{2}\times \left\{ 0\right\} \times $ $\mathbb{%
R}$ of $E$ and the group $A=\mathbb{R}^{2}\times $ $\mathbb{R}$ , direct
product of the group $\mathbb{R}^{2}$ by the group $\mathbb{R}$ with the
closed sub$-$group $\mathbb{R}^{2}\times \mathbb{R}\times \left\{ 0\right\} $%
\ of $E$

\bigskip \textbf{Definition 7.1. }\textit{For every\ }$\phi \in C^{\infty
}(G),$\textit{\ one can define a functions\ }$\tau \phi $ \textit{belong to }%
$C^{\infty }(G\times \mathbb{R}),$ \textit{and }$\iota \phi $ \textit{belong
to} $E$\textit{\ as follows:}%
\begin{equation}
\tau \phi (X,x,y)=\phi (x^{-1}X,x+y)
\end{equation}%
\begin{equation}
\iota \phi (X,x,y)=\phi (xX,x+y)
\end{equation}%
\textit{for any\ }$(X,x,y)\in G\times \mathbb{R}^{m}.\bigskip $ \textit{The
functions}\ $\tau \phi $ \textit{and }$\iota \phi $\ \textit{are invariant
in the following sense}%
\begin{equation}
\tau \phi (kX,x+k,y-k)=\phi (z,x,y)
\end{equation}%
\begin{equation}
\iota \phi (kX,x-k,y+k)=\phi (z,x,y)
\end{equation}

\bigskip Now, I state my theorem

\textbf{Theorem 7.1. }\textit{Let }$P$ \textit{be a right invariant
differential on} $G,$ \textit{and let} $u$ \textit{be the} \textit{%
distribution associated to} $P.$ \textit{Then the equation} 
\begin{equation}
P\phi (X,x)=u\ast \phi (X,x)=\int_{G}\phi ((w,v)^{-1}(X,x)u(X,x)dwdv=\varphi
(X,x)
\end{equation}%
\textit{has a solution} $\phi \in C^{\infty }(G),$\textit{for any function }$%
\varphi \in $ $C^{\infty }(G),$\textit{\ where }$\ast $ \textit{signifies
the convolution product on }$G.$

\textit{Proof: }Consider the operator $P$ as a differential operator $Q$ on
the abelian group $A=\mathbb{R}^{2}\times \left\{ 0\right\} \times \mathbb{R}%
.$ By the theory of partial differential equations with constant
coefficients on $\mathbb{R}^{2}\times \mathbb{R},$ then for any function $%
g\in C^{\infty }(\mathbb{R}^{2}\times \mathbb{R}),$\ there exist a function $%
\psi $ on $\mathbb{R}^{2}\times \mathbb{R},$ such that%
\begin{equation}
Q\psi (X,x)=u\ast _{c}\psi (X,x)=\int_{\mathbb{R}^{3}}\psi
(X-a,x-b)u(X,x)dadb=g(X,x)
\end{equation}

Using the extension of the function $\psi $ on the group $G\times \mathbb{R}%
, $then for each $f\in C^{\infty }(\mathbb{R}^{2}\times \mathbb{R}),$ I get 
\begin{equation}
=(u\ast _{c}\tau \psi )(X,0,y)\downarrow _{A}=f(X,y)
\end{equation}

\bigskip Let $\tau f$ be the extension of the function $f$ \ on the group $%
G\times \mathbb{R},$ that means 
\begin{eqnarray}
Q\tau \psi (X,0,y) &\downarrow &_{A}=(u\ast _{c}\tau \psi )(X,0,y)\downarrow
_{A}= \\
\tau f(X,0,y) &\downarrow &_{A}=f(X,y)  \notag
\end{eqnarray}%
where 
\begin{eqnarray}
(u\ast _{c}\psi )(X,y) &=&\int_{\mathbb{R}^{3}}\psi (X-a,y-b)u(X,y)dadb \\
&=&\tau f(X,0,y)\downarrow _{A}=f(X,y)
\end{eqnarray}

Let $\top _{x}$ be the right translation of the group $G$, which is defined
as%
\begin{equation}
=\top _{x}\Psi (X,t)=\Psi ((X,t)((0,x))=\Psi (X,t+x)
\end{equation}

Then I have $\tau \psi $ is the solution of the equation%
\begin{equation}
(u\ast \tau \psi )(X,x,0)\downarrow _{G}=f(x^{-1}X,x)
\end{equation}

In fact, we have 
\begin{eqnarray}
&=&\top _{x}(u\ast _{c}\tau \psi )(X,0,0)  \notag \\
&=&(u\ast \tau \psi )(X,x,0)\downarrow _{G}=(u\ast \tau \psi )(X,x,0)  \notag
\\
&=&\top _{x}\tau f(X,0,0)\downarrow _{G}=\tau f(X,x,0)=f(x^{-1}X,x)
\end{eqnarray}

So I get, if $\psi $ is the solution of the equation on the abelian group $A=%
\mathbb{R}^{2}\times \mathbb{R}$ 
\begin{equation}
(Q\psi )(X,y)=f(X,y)
\end{equation}%
on the abelian group $A=\mathbb{R}^{2}\times \mathbb{R}$, then the function $%
\tau \psi $ is the solution of the equation 
\begin{equation}
(P\tau \psi )(X,x)=f(x^{-1}X,x)
\end{equation}%
on the group $G.$ Let$\ \widetilde{\psi }(X,x)$ be the function, which is
defined as%
\begin{equation}
\widetilde{\psi }(X,x)=\psi (xX,x)
\end{equation}

In the same way, I have proved by in $[6],$ if $\widetilde{\psi }(X,x)$ is
the solution of the equation 
\begin{equation}
Q\text{ }\widetilde{\psi }(X,x)=\widetilde{\varphi }(X,x)
\end{equation}%
on the group $A,$ then the function $\psi $ is the solution of the equation%
\begin{equation}
P\psi (X,x)=\varphi (X,x)
\end{equation}%
on the group $G.$

\textbf{Corollary 7.1. }\textit{The Lewy equation is solvable in the sense,
for any }$g\in C^{\infty }(\mathbb{R}^{3})$ \textit{there is a function }$%
f\in C^{\infty }(\mathbb{R}^{3}),$ \textit{such that }%
\begin{equation}
L=(-\frac{\partial }{\partial x}-i\frac{\partial }{\partial y}+2i\text{ }(x+i%
\text{ }y)\frac{\partial }{\partial z})f=g
\end{equation}

The Lewy equation is invariant on the $3-$ dimensional nilpotent symplectic
group $N=\mathbb{R}^{2}\rtimes _{\rho }\mathbb{R}$. So it is solvable.

\textbf{The Example of Hormander for the non solvability}

Homander had considered in his book $[16,p.156],$ another form of the Lewy
operator, which is\textbf{\ }%
\begin{equation}
P(x,D)=(-i\partial _{x}+\partial _{y}-2x\partial _{z}-2iy\partial _{z})
\end{equation}

He constructed his example the operator of real variable coefficients, which
is%
\begin{equation}
Q(x,D)=P(x,D)\overline{P(x,D)}\overline{P(x,D)}P(x,D)
\end{equation}%
and proved $Q(x,D)$ is unsolvable see $[16,p164]$, where $\overline{P(x,D)}$
is the operator defined by 
\begin{equation}
\overline{P(x,D)}=(i\partial _{x}+\partial _{y}-2x\partial _{z}+2iy\partial
_{z})
\end{equation}

My result is:

\textbf{Theorem 7.2. }\textit{The operator }$Q(x,D)$ \textit{is solvable}

\textit{Proof: }Let $R$ be the following Cauchy-Riemann operator 
\begin{equation}
R=-i\partial _{x}+\partial _{y}
\end{equation}%
and let $\phi $\ be any function infinitely differentiable on $\mathbb{R}%
^{3} $, then we get

\begin{eqnarray}
&=&\hbar (-i\partial _{x})\hbar \phi \ (z,y,-x)=-i\partial _{x}\hbar \varphi
(z+2yx,y,x)  \notag \\
&=&(-i\frac{d}{dt})_{0}\hbar \phi \ (z+2yx,y,x+t)  \notag \\
&=&(-i\frac{d}{dt})_{0}\phi \ (z-2yt),y,-x-t)  \notag \\
&=&(-i\partial _{x}-2yi\partial _{z})\phi \ (z,y,-x)
\end{eqnarray}%
and%
\begin{eqnarray*}
\hbar (\partial _{y})\hbar \phi \ (z,y,-x) &=&\partial _{y}\hbar \phi \
(z+2xy,y,x) \\
&=&(\frac{d}{ds})_{0}\hbar \phi \ (z+2yx,y+s,x) \\
&=&(\frac{d}{ds})_{0}\phi \ (z-2sx,y+s,-x) \\
&=&(\partial _{y}-2x\partial _{z})\phi \ (z,y,-x)
\end{eqnarray*}

So, we get

\begin{equation}
(P(x,D)\phi )(z,y,-x)=(-i\partial _{x}+\partial _{y}-2x\partial
_{z}-2iy\partial _{z})\phi =\hbar R\hbar \phi (z,y,-x)\ 
\end{equation}

In the same manner, I prove

\begin{equation}
(\overline{P(x,D)}\phi )(z,y,-x)=(i\partial _{x}+\partial _{y}-2x\partial
_{z}+2iy\partial _{z})\phi =\hbar R_{\star }\hbar \phi (z,y,-x)\ 
\end{equation}%
where $R_{\star }$%
\begin{equation}
R_{\star }=i\partial _{x}+\partial _{y}
\end{equation}

Finally, I find 
\begin{equation}
((\overline{P(x,D)}(P(x,D)\phi )(z,y,-x)=\hbar R_{\star }R\hbar \phi
(z,y,-x)\ 
\end{equation}

\begin{equation}
((P(x,D)\overline{(P(x,D)})\phi )(z,y,-x)=\hbar RR\star \hbar \phi (z,y,-x)\ 
\end{equation}

\begin{eqnarray}
&&((((P(x,D)\overline{P(x,D)}\overline{P(x,D)}))))(P(x,D)\phi )(z,y,-x) 
\notag \\
&=&\hbar RR_{\star }R_{\star }R\hbar \phi (z,y,-x)\ =Q(x,D)\phi (z,y,-x)
\end{eqnarray}

Hence the solvability of the operator $Q(x,D).$ Also the operator 
\begin{equation}
X+iY-4iZ=ix\frac{\partial }{\partial _{z}}+i\frac{\partial }{\partial y}-y%
\frac{\partial }{\partial _{z}}+\frac{\partial }{\partial x}-4i\frac{%
\partial }{\partial _{z}}
\end{equation}%
is solvable. So the invalidity of the Hormander condition for the non
solvability

\section{\protect\bigskip \textbf{Conclusion}}

\textbf{8.1.} Any invariant differential operator has the form%
\begin{equation}
P=\sum\limits_{\alpha ,\beta }a_{\alpha ,\beta }X^{\alpha }Y^{\beta }
\end{equation}%
on the Lie group $G=\mathbb{R}^{2}\times _{\rho }\mathbb{R}$, where $%
X^{\alpha }=(X_{1}^{\alpha _{1}},X_{2}^{\alpha _{2}}),Y^{\beta },$ $\alpha
_{i}\in 
\mathbb{N}
\in 
\mathbb{N}
$ $(1\leq i\leq 2)$ and $X=(X_{1},X_{2}),$ are the invariant vectors field
on $G,$ which are the basis of the Lie algebra \underline{$g$} of $G$ and $%
a_{\alpha ,\beta }\in 
\mathbb{C}
.$ Any invariant partial differential equation on the $3-$dimensional group $%
G=\mathbb{R}^{2}\times _{\rho }\mathbb{R}$ is solvable. So the invalidity of
the Hormander condition for the non existence.

Over fifty years ago where there are a lot of books and lot of published
papers by many mathematicians as $[2,5,9,17,18,20,21,25],$are all based on a
non careful mathematical ideas. Especially those research published after $%
2006$ the date of opening my new way in Fourier analysis on non abelian Lie
groups. Unfortunately, some of those research books and articles were
published in the famous scientific centers such as Springer $[3,16]$,
Elsevier $[23]$, AMS $[22,27]$, Wiley $[26]$, Francis \& Taylor $[1]$, ...
ect

\textbf{8.2. Open questions. }The operator $[3,$ $p.2]$ 
\begin{eqnarray}
&&(y^{2}-z^{2})\frac{\partial ^{2}u}{\partial x^{2}}+(1+x^{2})(\frac{%
\partial ^{2}u}{\partial y^{2}}-\frac{\partial ^{2}u}{\partial z^{2}})-xy%
\frac{\partial ^{2}u}{\partial x\partial y}-  \notag \\
&&-\frac{\partial ^{2}(xyu)}{\partial x\partial y}+xz\frac{\partial ^{2}u}{%
\partial x\partial z}+\frac{\partial ^{2}(xyu)}{\partial x\partial z}
\end{eqnarray}%
can be solved

\bigskip \textbf{8.3. Open questions. }Is the operator

\begin{equation}
X^{2}+Y^{2}-4iZ=(x\frac{\partial }{\partial _{z}}+\frac{\partial }{\partial y%
})^{2}+(-y\frac{\partial }{\partial _{z}}+\frac{\partial }{\partial x}%
)^{2}-4i\frac{\partial }{\partial _{z}}
\end{equation}
solvable and hypoelliptic, and the operator

\begin{equation}
L-\alpha i\partial _{z}=\partial _{x}+2y\partial _{z}+i\partial
_{y}-2ix\partial _{z}-\alpha i\partial _{z},\alpha \in \mathbb{R}
\end{equation}%
is hypoelliptic

\textbf{Open} \textbf{Question. }Consider the Kannai operators $[3,$ $p.5]$%
\textbf{\ }%
\begin{equation}
D_{1}=\frac{\partial }{\partial _{x}}+x\frac{\partial ^{2}}{\partial _{y^{2}}%
},\text{ }D_{2}=\frac{\partial }{\partial _{x}}-x\frac{\partial ^{2}}{%
\partial _{y^{2}}}
\end{equation}%
I believe, the first can be solved on the $3-$dimensional Heisenberg group $%
H $ and the second can be hypoelliptic.

\end{document}